\documentclass[12pt]{article}
\usepackage{amsfonts,amssymb}

\def\C{\mathbb{C}}

\def\Z{\mathbb{Z}}

\def\g{\mathfrak g}

\def\bq{ \begin{equation} }
\def\eq{ \end{equation} }
\def\ben{ \begin{eqnarray} }
\def\en{ \end{eqnarray} }

\def\frac#1#2{{#1\over #2}}

\def\on#1#2{\mathop{\vbox{\ialign{##\crcr\noalign{\kern2pt}
$\scriptstyle{#2}$\crcr\noalign{\kern2pt\nointerlineskip}
\kern-2pt$\hfil\displaystyle{#1}\hfil$\crcr}}}\limits}

\begin{document}

\baselineskip=15pt
\vspace{1cm} \noindent {\LARGE \textbf {Pairs of compatible
associative algebras, classical Yang-Baxter equation and quiver
representations}} \vskip1cm \hfill
\begin{minipage}{13.5cm}
\baselineskip=15pt {\bf
 Alexander Odesskii ${}^{1,} {}^{2}$ and
 Vladimir Sokolov ${}^{1}$} \\ [2ex]
{\footnotesize ${}^1$ Landau Institute for Theoretical Physics,
Moscow, Russia
\\
${}^{2}$  School of Mathematics, The University of
Manchester, UK}\\
\vskip1cm{\bf Abstract}

Given an associative multiplication in matrix algebra compatible
with the usual one or, in other words, linear deformation of matrix
algebra, we  construct a solution to the classical Yang-Baxter
equation. We also develop a theory of such deformations and
construct numerous examples. It turns out that these deformations
are in one-to-one correspondence with representations of certain
algebraic structures, which we call $M$-structures. We also describe
an important class of $M$-structures related to the affine Dynkin
diagrams of $A$, $D$, $E$-type. These $M$-structures and their
representations are described in terms of quiver representations.

\end{minipage}

\vglue1cm \textbf{Address}: Landau Institute for Theoretical
Physics, Kosygina 2, 119334, Moscow, Russia

\textbf{E-mail}: odesskii@itp.ac.ru, \, sokolov@itp.ac.ru \newpage

\centerline{\Large\bf Introduction}
\medskip
 Two
associative algebras with multiplications $(a,b)\to ab$ and
$(a,b)\to a\circ b$ defined on the same finite dimensional vector
space are said to be compatible if the multiplication
\begin{equation} \label{pensil}
a \bullet b= a b+ \lambda \,a\circ b
\end{equation}
is associative for any constant $\lambda$. The multiplication
$\bullet$ can be regarded as a deformation of the multiplication
$(a,b)\to ab$ linear in the parameter $\lambda$.

In \cite{odsok1} we have studied multiplications  compatible with
the matrix product  or, in other words, linear deformations of the
matrix multiplication. It turns out that these deformations of the
matrix algebra are in one-to-one correspondence with representations
of certain algebraic structures, which we call $M$-structures. The
case of direct sum of several matrix algebras corresponds to
representations of the so-called $PM$-structures (see
\cite{odsok1}).

Given a pair of compatible associative products, one can construct a
hierarchy of integrable systems of ODEs via Lenard-Magri scheme
\cite{magri}. The Lax representations for these systems are
described in \cite{golsok4}. If one of the multiplications is the
usual matrix product, the integrable systems are Hamiltonian
$gl(N)$-models with quadratic Hamiltonians \cite{odsok2}. These
systems can be regarded as a generalization of the matrix equations
considered in \cite{sokmih}. Their skew-symmetric reductions give
rise to new integrable quadratic $so(n)$-Hamiltonians.

The main ingredient of the $M$-structure is a pair of associative
algebras $\cal A$ and $\cal B$ of the same dimension. The simplest
version of a structure of this kind can be regarded as an
associative analog of the Lie bi-algebra.

We define an {\it associative bi-algebra} as a couple of associative
algebras $\cal A$ and $\cal B$ with non-degenerated pairing and a
structure of ${\cal B}\otimes{\cal A}^{op}$-module on the space
$\cal L=\cal A\oplus\cal B$ such that the algebra $\cal A$ acts on
$\cal A\subset \cal L$ by right multiplications, the algebra $\cal
B$ acts on $\cal B\subset \cal L$ by left multiplications and the
pairing is invariant with respect to this action (that is
$(bb^{\prime},a)=(b,b^{\prime}a)$ and
$(b,aa^{\prime})=(ba,a^{\prime})$ for $a,a^{\prime}\in\cal A$ and
$b,b^{\prime}\in\cal B$). Here ${\cal A}^{op}$ stands for the
algebra opposite to $\cal A$. Given an associative bi-algebra one
has the structure of associative algebra in the space $\cal
A\oplus\cal B\oplus\cal A\otimes\cal B$ (this is an analog of the
Drinfeld double).

In this paper we introduce the notion of {\it associative
$r$-matrices}, which is a particular case of usual classical
$r$-matrices. It turns out that the constant associative
$r$-matrices can be classified in terms of associative bi-algebras.
Moreover, one can introduce spectral parameters to the definition of
associative bi-algebras and obtain a classification of non-constant
associative $r$-matrices.

In \cite{odsok1} we have discovered an important class of $M$ and
$PM$-structures. These structures are related to the Cartan matrices
of affine Dynkin diagrams of the $\tilde A_{2 k-1},$ $\tilde D_{k},$
$\tilde E_{6},$ $\tilde E_{7},$ and $\tilde E_{8}$-type. In this
paper we describe these $M$-structures and their representations in
terms of quiver representations.

The paper is organized as follows. In Section {\bf 1}, we consider
an associative analog of the classical Yang-Baxter equation. Since
semi-simple associative algebras are more rigid algebraic structures
then semi-simple Lie algebras, it turns out to be possible to
construct a developed theory of the associative Yang-Baxter equation
in the semi-simple case. This theory is suitable for constructing of
wide class of solutions to the Yang-Baxter equation. We are planning
to write a separate paper devoted to systematic search for
solutions.

In Section {\bf 2}, we give an explicit construction of a solution
to the Yang-Baxter equation by each pair of compatible Lie brackets
provided that the first bracket is rigid. The corresponding
$r$-matrices are not unitary and therefore they are out of
classification by A.~Belavin and V.~Drinfeld \cite{beldr}. In
particular, compatible associative products give rise to solutions
of the associative Yang-Baxter equation. This gives us a way to
construct $r$-matrices related to $M$-structures.

In Section {\bf 3} we recall the notion of $M$-structure and
formulate main results describing the relationship between
associative multiplications in matrix algebra compatible with the
usual matrix product and $M$-structures.

In  Section {\bf 4} we describe all $M$-structures with semi-simple
algebras $\cal A$ and $\cal B$. It turns out that such
$M$-structures are related to the Cartan matrices of affine Dynkin
diagrams of the $\tilde A_{2 k-1},$ $\tilde D_{k},$ $\tilde E_{6},$
$\tilde E_{7},$ and $\tilde E_{8}$-type. We describe these
$M$-structures and their representations in terms of representations
of affine quivers \cite{quivers1, quivers2, quivers3}.

In Appendix we give explicit formulas for these $M$-structures,
their representations and for corresponding solutions to the
classical Yang-Baxter equation.

\section{Classical Yang-Baxter equation}

Let $\g$ be a Lie algebra. Let $r(u,v)$ be a meromorphic function of
two complex variables with values in $End(\g)$. For each $u\in\C$ we
denote by $\g_u$ a vector space canonically isomorphic to $\g$. Let
$\tilde\g=\oplus_u\g_u$. We define a bracket on the space $\tilde\g$
by the formula
\begin{equation} \label{brinf}
[x_u,y_v]=([x,r(u,v)y])_u+([r(v,u)x,y])_v.
\end{equation}

{\bf Lemma 1.1.} The bracket (\ref{brinf}) defines a structure of a
Lie algebra on $\tilde\g$ iff $r(u,v)$ satisfies the following
equation
\begin{equation} \label{yblie}
[r(u,w)x,r(u,v)y]-r(u,v)[r(v,w)x,y]-r(u,w)[x,r(w,v)y]\in Cent(\g),
\end{equation}
where $x,y$ are arbitrary elements of $\g$ and  $Cent(\g)$ stands for the center of $\g$.

{\bf Proof} of the lemma is straightforward.

{\bf Definition.} The operator relation
\begin{equation} \label{yblie1}
[r(u,w)x,r(u,v)y]-r(u,v)[r(v,w)x,y]-r(u,w)[x,r(w,v)y]=0.
\end{equation}
is called {\it classical Yang-Baxter equation}. A solution $r(u,v)$
to the classical Yang-Baxter equation is called {\it classical
}$r$-{\it matrix}. Arguments of $r(u,v)$ are called spectral
parameters.

Note that the arguments $u,v$ of $r$ could be also elements of
$\C^n$ for $n>1$ or elements of some complex manifold, which is
called the manifold of spectral parameters.

Suppose $\g$ possesses a non-degenerate invariant scalar product
$(\cdot,\cdot)$. Then it is easy to see \cite{sem} that the equation
 (\ref{yblie1}) is equivalent to the classical Yang-Baxter
equation written in the tensor form. An $r$-matrix is called unitary
if $(x,r(u,v)y)=-(r(v,u)x,y)$. The unitary $r$-matrices were
classified in \cite{beldr}. There is no any classification of
$r$-matrices in the general case.

It turns out that a theory of (non-unitary) $r$-matrices can be
developed in the special case of associative algebras. Let $A$ be an
associative algebra. Let $r(u,v)$ be a meromorphic function in two
complex variables with values in $End(A)$. For each $u\in\C$ we
denote by $A_u$ a vector space canonically isomorphic to $A$. Let
$\tilde A=\oplus_uA_u$. We define a product on the space $\tilde A$
by the formula
\begin{equation} \label{assinf}
x_uy_v=(x(r(u,v)y))_u+((r(v,u)x)y)_v.
\end{equation}

{\bf Lemma 1.2.} The product (\ref{assinf}) defines a structure of
an associative algebra on $\tilde A$ iff $r(u,v)$ satisfies the
following equation
\begin{equation} \label{ybass}
(r(u,w)x)(r(u,v)y)-r(u,v)((r(v,w)x)y)-r(u,w)(x(r(w,v)y))\in Null(A),
\end{equation}
where $Null(A)$ is the set of $z\in A$ such that $zt=tz=0$ for all
$t\in A$.

{\bf Proof} of the lemma is straightforward.

{\bf Definition.} The relation
\begin{equation} \label{ybass1}
(r(u,w)x)(r(u,v)y)-r(u,v)((r(v,w)x)y)-r(u,w)(x(r(w,v)y))=0
\end{equation}
is called {\it associative Yang-Baxter equation}.

{\bf Lemma 1.3.} Let $\g$ be a Lie algebra with the brackets
$[x,y]=xy-yx$. Then any solution of (\ref{ybass1}) is a solution of
(\ref{yblie1}).

{\bf Proof} of the lemma is straightforward.

Let $A=Mat_n$. It is easy to see that any operator from $End(A)$ to
$End(A)$ has the form $x\to a_1\,x\,b^1+...+a_p\,x\,b^p$ for some
matrices $a_1,...,a_p,b^1,...,b^p$. Moreover, $p$ is the smallest
possible for such representation iff the matrices $\{a_1,...,a_p\}$
are linear independent as well as $\{b_1,...,b_p\}$.

{\bf Theorem 1.1.} Let
$$r(u,v)x=a_1(u,v)\,x\,b^1(v,u)+...+a_p(u,v)\,x\,b^p(v,u),$$ where
$a_1(u,v),...,b^p(u,v)$ are meromorphic functions with values in
$Mat_n$ such that \linebreak $\{a_1(u,v),...,a_p(u,v)\}$ are linear
independent over the field of meromorphic functions in $u,v$ as well
as $\{b^1(u,v),...,b^p(u,v)\}$. Then $r(u,v)$ satisfies
(\ref{ybass1}) iff there exist meromorphic functions
$\phi_{i,j}^k(u,v,w)$ and $\psi_{i,j}^k(u,v,w)$ such that
$$a_i(u,v)a_j(v,w)=\phi_{i,j}^k(u,v,w)a_k(u,w),$$
\begin{equation} \label{bialg}
b^i(u,v)b^j(v,w)=\psi^{i,j}_k(u,v,w)b^k(u,w),
\end{equation}
$$b^i(u,v)a_j(v,w)=\phi^i_{j,k}(v,w,u)b^k(u,w)+\psi_j^{k,i}(w,u,v)a_k(u,w).$$
The tensors $\phi_{i,j}^k(u,v,w)$ and $\psi_{i,j}^k(u,v,w)$ satisfy
the following equations
$$\phi_{i,j}^s(u,v,w)\phi_{s,k}^l(u,w,t)=\phi_{i,s}^l(u,v,t)\phi_{j,k}^s(v,w,t),$$
\begin{equation} \label{eq}
\psi^{i,j}_s(u,v,w)\psi^{s,k}_l(u,w,t)=\psi^{i,s}_l(u,v,t)\psi^{j,k}_s(v,w,t),
\end{equation}
$$\phi_{j,k}^s(v,w,t)\psi_s^{l,i}(t,u,v)=
\phi_{s,k}^l(u,w,t)\psi_j^{s,i}(w,u,v)+\phi_{j,s}^i(v,w,u)\psi_k^{l,s}(t,u,w).$$

{\bf Proof} of the theorem is similar to the proof of the Theorem
3.1 from \cite{odsok1}.

{\bf Remark 1.} It is easy to give an invariant description of the
corresponding algebraic structure. In the case of constant
$r$-matrix this leads to the associative bi-algebras described in
the Introduction.

{\bf Remark 2.} A similar statement holds in the case of semi-simple
algebra $A$.

{\bf Example 1.} Let $r(u,v) x=\frac{1}{u-v}e(u,v) x f(v,u),$ where
$$e(u,v)e(v,w)=e(u,w),\qquad f(u,v)f(v,w)=f(u,w),$$
\begin{equation} \label{eqef}
e(u,v)f(v,w)=\frac{u-v}{u-w}f(u,w)+\frac{v-w}{u-w}e(u,w).
\end{equation}
These equations hold if we assume, for example, that $e(u,v)=1,\,\,
f(u,v)=(v-u)(u+C)^{-1}+1,$ where $C$ is an arbitrary constant
matrix.

{\bf Example 2.} Let $A=\C^p$. The algebra $A$ has a basis
$\{e_i,i=1...p\}$ such that $e_ie_j=\delta_{i,j}e_i$. The formula
$$r(u,v)e_i=\sum_{1\leq j\leq
p}\frac{\psi_i(v)}{\phi_j(u)-\phi_i(v)}e_j$$ gives an associative
$r$-matrix for any functions $\phi_1,...,\phi_p,\psi_1,...,\psi_p$
in one variable, where $\phi_1,...,\phi_p$ are not constant. This
$r$-matrix can be written in the form
$$r(\vec{u},\vec{v})=\sum_{1\leq j\leq
p}\frac{\psi_i(\vec{v})}{u_j-v_i}e_j,$$ where
$\vec{u}=(u_1,...,u_p)$, $\vec{v}=(v_1,...,v_p)$, $\psi_i(\vec{v})$
are functions in $p$ variables. In this case, the manifold of
spectral parameters is $\C^p$.

\section{Compatible products and solutions to the classical Yang-Baxter equation}

Two Lie brackets $[\cdot,\cdot]$ and $[\cdot,\cdot]_1$ defined on
the same vector space  $\g$ are said to be compatible if
$[\cdot,\cdot]_{\lambda}=[\cdot,\cdot]+\lambda [\cdot,\cdot]_1$ is a
Lie bracket for any $\lambda$. In the papers
\cite{golsok1,golsok2,golsok3,odsok} different applications of the
notion of compatible Lie brackets to the integrability theory have
been considered.

Suppose that the bracket $[\cdot,\cdot]$ is rigid, i.e. $H^2(\g
,\g)=0$ with respect to $[\cdot,\cdot]$. In this case the Lie
algebras with the brackets $[\cdot,\cdot]_{\lambda}$ are isomorphic
to the Lie algebra with the bracket $[\cdot,\cdot]$ for almost all
values of the parameter $\lambda$. This means that there exists a
meromorphic function $\lambda\to S_{\lambda}$ with values in
$End(\g)$ such that $S_0=Id$ and
\begin{equation}\label{isomlie}
[S_{\lambda}(x),S_{\lambda}(y)]=S_{\lambda}([x,y]+\lambda [x, y]_1).
\end{equation}

{\bf Theorem 2.1.} The formula
\begin{equation}\label{ybsol}
r(u,v)=\frac{1}{u-v}S_uS_v^{-1}
\end{equation}
defines a solution to the classical Yang-Baxter equation
(\ref{yblie1}).

{\bf Proof.} For $r(u,v)$ given by (\ref{ybsol}) equation
(\ref{yblie1}) is equivalent to
\begin{equation} \label{ybsubs}
\frac{1}{(u-v)(u-w)}[S_uS_w^{-1}(x),S_uS_v^{-1}(y)]-\frac{1}{(u-v)(v-w)}S_uS_v^{-1}([S_vS_w^{-1}(x),y])
-
\end{equation}
$$\frac{1}{(u-w)(w-v)}S_uS_w^{-1}([x,S_wS_v^{-1}(y)])=0.$$
Using (\ref{isomlie}), we get
$$[S_uS_w^{-1}(x),S_uS_v^{-1}(y)]=S_u([S_w^{-1}(x),S_v^{-1}(y)]+u[S_w^{-1}(x), S_v^{-1}(y)]_1),$$
$$S_uS_v^{-1}([S_vS_w^{-1}(x),y])=S_u([S_w^{-1}(x),S_v^{-1}(y)]+v[S_w^{-1}(x), S_v^{-1}(y)]_1),$$
$$S_uS_w^{-1}([x,S_wS_v^{-1}(y)])=S_u([S_w^{-1}(x),S_v^{-1}(y)]+w[S_w^{-1}(x), S_v^{-1}(y)]_1).$$
Substituting these expressions into the left hand side of
(\ref{ybsubs}), we obtain the statement.

{\bf Remark 1.} It is clear that the $r$-matrix (\ref{ybsol}) is
unitary with respect to an invariant form $(\cdot,\cdot)$ if the
operator $S_{\lambda}$ is orthogonal. In this case the formula
(\ref{isomlie}) implies that the form $(\cdot,\cdot)$ is invariant
with respect to the second bracket.

Two associative algebras with multiplications $(x,y)\to xy$ and
$(x,y)\to x\circ y$ defined on the same finite dimensional vector
space $A$ are said to be {\it compatible} if the multiplication
(\ref{pensil}) is associative for any constant $\lambda$.  Suppose
$H^2(A,A)=0$ with respect to the first multiplication; then there
exists a meromorphic function $\lambda\to S_{\lambda}$ with values
in $End(A)$ such that $S_0=Id$ and
\begin{equation}\label{isom}
S_{\lambda}(x)S_{\lambda}(y)=S_{\lambda}(xy+\lambda x\circ y).
\end{equation}
The Taylor decomposition of $S_{\lambda}$ at $\lambda=0$ has the
following form
\begin{equation}
S_{\lambda}=1+ R \ \lambda+ S \ \lambda^2 + \cdots \,.  \label{RS}
\end{equation}
Substituting this decomposition into (\ref{isom}) and equating the
coefficients of $\lambda,$ we obtain the formula
\begin{equation} \label{mult2}
x \circ y =R(x)y+xR(y)-R(xy),
\end{equation}
where $R$ is a linear operator on $A$. It is clear that  for any
$a\in A$ the transformation
\begin{equation}
     R\longrightarrow R + ad_a,
\label{ad}
\end{equation}
where $ad_a$ is a linear operator $v\to av-va,$ does not change the
multiplication $\circ.$

{\bf Definition.} Operators $R$ and $R^{\prime}$ are said to be
equivalent if $R-R^{\prime}=ad_a$ for some $a\in A$.

The following analog of the Theorem 2.1 can be proved similarly.

{\bf Theorem 2.2.} Suppose that $S_{\lambda}$ satisfies
(\ref{isom}), then the formula (\ref{ybsol}) defines a solution to
the associative Yang-Baxter equation (\ref{ybass1}).

{\bf Remark 2.} In the important particular case
$S_{\lambda}=1+\lambda R,$ the $r$-matrix (\ref{ybsol}) is
equivalent to
\begin{equation}\label{ybsol1}
r(u,v)=\frac{1}{u-v}+(v+R)^{-1}.
\end{equation}

 Let $A=Mat_N$. Consider the following classification problem: to describe all
possible associative multiplications $\circ$  compatible with the
usual matrix product in $A$. Since $H^2(A,A)=0$ for any semi-simple
associative algebra, an operator valued meromorphic function
$S_{\lambda}$ with the properties $S_0=Id$ and (\ref{isom}) exists
for any such multiplication and the multiplication is given by the
formula (\ref{mult2}).

{\bf Example.} Let $a\in Mat_N$ be an arbitrary matrix and $R$ be
the operator of left multiplication by $a$. Then (\ref{mult2})
yields the multiplication $x\circ y=x a y,$ which is associative and
compatible with the standard one. It is clear that $S_{\lambda}$ can
be chosen in the form $S_{\lambda}(x)=(1+\lambda a)x$. We have
$$r(u,v)=\frac{1}{u-v}+(v+ a)^{-1}$$
in this case.

Any linear operator $R$ on the space $Mat_N$ may be written in the
form $R(x)=a_1 x b^1+...+a_l x b^l$ for some matrices
$a_1,...,a_l,b^1,...,b^l$. Indeed, the operators $x\to e_{i,j} x
e_{i_1,j_1}$ form a basis in the space of linear operators on
$Mat_N$.

It is convenient to represent the operator $R$ from the formula
(\ref{mult2}) in the form
\begin{equation}\label{Rmat}
R(x)=a_1 \,x \,b^1+...+a_p \,x\, b^p+c\, x
\end{equation}
with $p$  being smallest possible in the class of equivalence of
$R$. This means that the matrices $\{a_1,...,a_p,1\}$ are linear
independent as well as the matrices $\{b^1,...,b^p,1\}$. According
to (\ref{mult2}), the second product has the following form
\begin{equation}\label{sekprod}
x\circ y=\sum_i (a_i\, x\, b^i\, y+x \,a_i\, y\, b^i-a_i\, x y\,
b^i)+x\, c\, y.
\end{equation}

It turns out that the matrices $\{a_1,...,a_p,b^1,...,b^p,c\}$ form
a representation of a certain algebraic structure. We describe this
structure in the next section.

\section{$M$-structures and the corresponding associative algebras}

In this Section we formulate the results of \cite{odsok1} and their
simple consequences we will use below.

{\bf Definition.} By weak $M$-structure on a linear space $\cal L$
we mean the following data:
\begin{itemize}
\item  Two subspaces $\cal A$ and $\cal B$ and a distinguished
element $1\in\cal A\cap\cal B\subset \cal L$.

\item  A non-degenerate symmetric scalar product $(\cdot, \cdot)$ on the
space $\cal L$.

\item  Associative products $\cal A\times\cal A\to\cal A$ and
$\cal B\times\cal B\to\cal B$ with unity $1$.

\item  A left action $\cal B\times\cal L\to\cal L$ of the algebra
$\cal B$ and a right action $\cal L\times\cal A\to\cal L$ of the
algebra $\cal A$ on the space $\cal L$ that commute to each other.
\end{itemize}

These data should satisfy the following properties:

{\bf 1.} $\dim{\cal A\cap\cal B=\dim\cal L/(\cal A+\cal B)}= 1$.

{\bf 2.} The restriction of the action $\cal B\times\cal L\to\cal L$
to the subspace $\cal B\subset \cal L$ is the product in $\cal B$.
The restriction of the action $\cal L\times\cal A\to\cal L$ to the
subspace $\cal A\subset \cal L$ is the product in $\cal A$.

{\bf 3.} $(a_1,a_2)=(b_1,b_2)=0$ and $(b_1b_2,v)=(b_1,b_2v)$,
$(v,a_1a_2)=(va_1,a_2)$ for any $a_1,a_2\in\cal A$, $b_1,b_2\in\cal
B$ and $v\in\cal L$.

It follows from these properties that $(\cdot,\cdot)$ defines a non
- degenerate pairing between ${\cal A}/\C1$ and ${\cal B}/\C1$ and
therefore $\dim\cal A=\dim\cal B$ and $\dim{\cal L}=2\dim{\cal A}$.

Given a weak $M$-structure $\cal L$ we define an associative algebra
$U(\cal L)$ generated by $\cal L$ and satisfying natural
compatibility and universality conditions.

{\bf Definition.} By weak $M$-algebra associated with a weak
$M$-structure $\cal L$ we mean an associative algebra $U(\cal L)$
with a linear mapping $j:{\cal L}\to U(\cal L)$ such that the
following conditions are satisfied:

{\bf 1.} $j(b)j(x)=j(bx)$ and $j(x)j(a)=j(xa)$ for $a\in\cal A$,
$b\in\cal B$ and $x\in\cal L$.

{\bf 2.} For any algebra $X$ with a linear mapping $j^{\prime}:{\cal
L}\to X$ satisfying the property {\bf 1} there exists a unique
homomorphism of algebras $f: U({\cal L})\to X$ such that $f\circ
j=j^{\prime}$.

It is easy to see that $U(\cal L)$ exists and is unique for given
$\cal L$.

{\bf Definition.} A weak $M$-structure $\cal L$ is called
$M$-structure if there exists a central element $K\in U(\cal L)$ of
the algebra $U(\cal L)$ quadratic with respect to $\cal L$.

{\bf Theorem 3.1.} Let $\cal L$ be an $M$-structure. Then there
exists a basis $\{1,A_1,...,A_p,B_1,..., B_p,C\}$ in $\cal L$ such
that $\{1,A_1,...,A_p\}$ is a basis in $\cal A$, $\{1,B_1,...,B_p\}$
is a basis in $\cal B$ and
$$K=A_1B_1+...+A_pB_p+C.$$

{\bf Theorem 3.2.} Let $R\in End(U(\cal L))$ is given by the
formula
$$R(x)=A_1xB_1+...+A_pxB_p+Cx$$ and $\circ$ is defined by (\ref{mult2}).
Then $\circ$ is associative and compatible with the usual product in
$U(\cal L)$.

Notice that $K=R(1)$.

{\bf Theorem 3.3.} Let $\circ$ be an associative product in the
space $Mat_N$ compatible with the usual one and written in the form
(\ref{mult2}), where $R$ is given by (\ref{Rmat}) with $p$ being
smallest possible in the class of equivalence of $R$. Then there
exists an $M$-structure $\cal L$ with representation $U({\cal L})\to
Mat_N$ such that $\dim{\cal A}=\dim{\cal B}=p+1$, image of $\cal A$
has the basis $\{1,a_1,...,a_p\}$ and image of $\cal B$ has the
basis $\{1,b_1,...,b_p\}$.

{\bf Definition.} A representation of  $U(\cal L)$ is called
non-degenerate if its restrictions on the algebras $\cal A$ and
$\cal B$ are exact.

{\bf Theorem 3.4.} There is one-to-one correspondence between $N$-
dimensional non-degenerate representations of algebras $U(\cal L)$
corresponding to $M$-structures and associative products in $Mat_N$
compatible with the usual matrix product.

The structure of the algebra $U(\cal L)$ for $M$-structure $\cal L$
can be described as follows.

{\bf Theorem 3.5.} The algebra $U(\cal L)$ is spanned by the
elements of the form $a\,b\,K^s,$ where $ a\in {\cal A}, b\in {\cal
B},\, s\in \Z_{+}$.

We need also the following

{\bf Definition.} Let ${\cal L}$ be a weak $M$-structure. By the
opposite weak $M$-structure ${\cal L}^{op}$ we mean an $M$-structure
with the same linear space ${\cal L}$, the same scalar product and
algebras $\cal A$, $\cal B$ replaced by the opposite algebras ${\cal
B}^{op}$, ${\cal A}^{op},$ correspondingly.

It is easy to see that if ${\cal L}$ is an $M$-structure, then
${\cal L}^{op}$ is an $M$-structure as well.

\section{$M$-structures with semi-simple algebras $\cal A$ and $\cal B$ and quiver representations}

\subsection{Matrix of multiplicities}

By $V^l$  we denote the direct sum of $l$ copies of a linear space
$V.$ By definition, we put $V^0=\{0\}$. Recall \cite{alg} that any
semi-simple associative algebra over $\C$ has the form $\oplus_{1\le
i\le r}End(V_i),$ any left $End(V)$-module has the form $V^l,$ and
any right $End(V)$-module has the form $(V^{\star})^l$ for some $r$
and $l$.

{\bf Lemma 4.1.} Let $\cal L$ be a weak $M$-structure. Suppose
${\cal A}=\oplus_{1\le i\le r}End(V_i),$ where $\dim V_i=m_i$. Then
$\cal L$ as a right $\cal A$-module is isomorphic to $\oplus_{1\le
i\le r}(V^{\star}_i)^{2m_i}$.

{\bf Proof.} Since any right $\cal A$-module has the form
$\oplus_{1\le i\le r}(V^{\star}_i)^{l_i}$ for some $l_1,...,l_r\ge
0,$ we have ${\cal L}=\oplus_{1\le i\le r}{\cal L}_i$ where ${\cal
L}_i=(V^{\star}_i)^{l_i}$. Note that $\cal A\subset\cal L$ and,
moreover, $End(V_i)\subset{\cal L}_i$ for $i=1,...,r$. Besides,
$End(V_i)\bot{\cal L}_j$ for $i\ne j$. Indeed, we have
$(v,a)=(v,Id_ia)=(vId_i,a)=0$ for $v\in{\cal L}_j$ and $a\in
End(V_i),$ where $Id_i$ is the unity of the subalgebra $End(V_i)$.
Since $(\cdot,\cdot)$ is non-degenerate and $End(V_i)\bot End(V_i)$
by the property 3 of weak $M$-structure, we have $\dim{\cal L}_i\ge
2\dim End(V_i)$. But $\sum_i\dim{\cal L}_i=\dim{\cal L}=2\dim{\cal
A}=\sum_i 2\dim End(V_i)$ and we obtain  $\dim{\cal L}_i=2\dim
End(V_i)$ for each $i=1,...,r,$ which is equivalent to the statement
of Lemma 4.1.

{\bf Lemma 4.2.} Let $\cal A$ and $\cal B$ be semi-simple:
\begin{equation}\label{AB}
{\cal A}=\oplus_{1\le i\le r}End(V_i),\qquad {\cal B}=\oplus_{1\le
j\le s}End(W_j), \qquad \dim V_i=m_i, \quad\dim W_j=n_j.
\end{equation}
 Then $\cal
L$ as an $\cal A\otimes\cal B$-module is given by the formula
\begin{equation}\label{L}
{\cal L}=\oplus_{1\le i\le r,1\le j\le s}(V^{\star}_i\otimes
W_j)^{a_{i,j}},
\end{equation}
where $a_{i,j}\ge 0$ and
\begin{equation}\label{adm}   \sum_{j=1}^s
a_{i,j}n_j=2m_i,\qquad \sum_{i=1}^r a_{i,j}m_i=2n_j.
\end{equation}

{\bf Proof.} It is known that any ${\cal A}^{op}\otimes\cal
B$-module has the form $\oplus_{1\le i\le r,1\le j\le
s}(V^{\star}_i\otimes W_j)^{a_{i,j}},$ where $a_{i,j}\ge 0$.
Applying Lemma 4.1, we obtain $\dim {\cal L}_i=2m_i^2,$ where ${\cal
L}_i=\oplus_{1\le j\le s}(V^{\star}_i\otimes W_j)^{a_{i,j}}$. This
gives the first equation from (\ref{adm}). The second equation can
be obtained similarly.

{\bf Definition.} The $r\times s$-matrix $(a_{i,j})$ from the Lemma
4.2 is called the matrix of multiplicities of the weak $M$-structure
$\cal L$.

{\bf Definition.} The $r\times s$-matrix $(a_{i,j})$ is called
 decomposable  if
there exist partitions $\{1,...,r\}=I\sqcup I^{\prime}$ and
$\{1,...,s\}=J\sqcup J^{\prime}$ such that $a_{i,j}=0$ for
$(i,j)\in I\times J^{\prime}\sqcup I^{\prime}\times J$.

{\bf Lemma 4.3.} The matrix of multiplicities is indecomposable.

{\bf Proof.} Suppose $(a_{i,j})$ is decomposable. We have $\cal
A=\cal A^{\prime}\oplus\cal A^{\prime\prime}$, $\cal B=\cal
B^{\prime}\oplus\cal B^{\prime\prime}$ and $\cal L=\cal
L^{\prime}\oplus\cal L^{\prime\prime},$ where $${\cal
A}^{\prime}=\oplus_{i\in I}End(V_i),\qquad {\cal
A}^{\prime\prime}=\oplus_{i\in I^{\prime}}End(V_i),\qquad {\cal
B}^{\prime}=\oplus_{j\in J}End(W_j),$$ $${\cal
B}^{\prime\prime}=\oplus_{j\in J^{\prime}}End(W_j), \qquad {\cal
L}^{\prime}=\oplus_{(i,j)\in I\times J}(V^{\star}_i\otimes
W_j)^{a_{i,j}},\qquad {\cal L}^{\prime\prime}=\oplus_{(i,j)\in
I^{\prime}\times J^{\prime}}(V^{\star}_i\otimes W_j)^{a_{i,j}}.$$
Let $1=e_1+e_2,$ where $e_1\in{\cal L}^{\prime}$ and $e_2\in{\cal
L}^{\prime\prime}$. It is clear that $e_1,e_2\in  \cal A\cap\cal B$.
Therefore, $\dim{\cal A\cap\cal B}>1,$ which contradicts to the
property 1 of weak $M$-structure.

Note that if $A$ is the matrix of multiplicities of a weak $M$
structure with semi-simple algebras $\cal A$ and $\cal B$, then
$A^t$ is the matrix of multiplicities of the opposite weak
$M$-structure.

{\bf Theorem 4.1.} Let $\cal L$ be a weak $M$-structure with
semi-simple algebras $\cal A$ and $\cal B$ given by the formula
(\ref{AB}) and with $\cal L$ given by (\ref{L}). Then there exists a
simple laced affine Dynkin diagram \cite{burb} with vector spaces
from the set $\{V_1,...,V_r,W_1,...,W_s\}$ assigned to each vertex
in such a way that:

{\bf 1.} there is an one-to-one correspondence between this set and
the set of vertices,

{\bf 2.} for any $i,j$ the spaces $V_i$, $V_j$ are not connected by
edges  as well as the spaces $W_i$, $W_j$,

{\bf 3.}  $a_{i,j}$ is equal to the number of edges between $V_i$
and $W_j$,

{\bf 4.} the vector $(\dim V_1,...,\dim V_r,\dim W_1,...,\dim W_s)$
is a positive imaginary root of the diagram.

{\bf Proof.} Consider a linear space with a basis
$\{v_1,...,v_r,w_1,...,w_s\}$ and the symmetric bilinear form
$(v_i,v_j)=(w_i,w_j)=2\delta_{i,j}$, $(v_i,w_j)=-a_{i,j}$. Let
$J=m_1 v_1+...+m_r v_r+n_1 w_1+...+n_s w_s$. It is clear that the
equations (\ref{adm}) can be written as $(v_i,J)=(w_j,J)=0,$ which
means that $J$ belongs to the kernel of the form $(\cdot,\cdot)$.
Therefore (see \cite{vinberg}), the matrix of the form is a Cartan
matrix of a simple laced affine Dynkin diagram. It is also clear
that $J$ is a positive imaginary root.

On the other hand, consider a simple laced affine Dynkin diagram
with a partition of the set of vertices into two subsets such that
vertices of the same subset are not connected. It is clear that if
such a partition exists, then it is unique up to transposition of
subsets. Let $v_1,...,v_r$ be roots corresponding to vertices of the
first subset and $w_1,...,w_s$ be roots corresponding to the second
subset. We have $(v_i,v_j)=(w_i,w_j)=2\delta_{i,j}$. Let
$a_{i,j}=-(v_i,w_j)$ and $J=m_1v_1+...+m_rv_r+n_1w_1+...+n_sw_s$ be
an imaginary root. It is clear that (\ref{adm}) holds.

{\bf Remark.} The interchanging  of the subsets corresponds to the
transposition of matrix $(a_{i,j})$.

It is easy to see that among simple laced affine Dynkin diagrams
only diagrams of the $\tilde A_{2 k-1},$ $\tilde D_{k},$ $\tilde
E_{6},$ $\tilde E_{7},$ and $\tilde E_{8}$-type admit a partition of
the set of vertices into two subsets such that vertices of the same
subset are not connected.  The natural question arises: to describe
all $M$-structures with the algebras $\cal A$ and $\cal B$ given by
(\ref{AB}) and $\cal L$ given by (\ref{L}), where the matrix
$(a_{i,j})$ is constructed by the affine Dynkin diagram of the
$\tilde A_{2 k-1},$ $\tilde D_{k},$ $\tilde E_{6},$ $\tilde E_{7},$
and $\tilde E_{8}$-type. It turns out that these $M$-structures
exist iff $J$ is the minimal positive imaginary root.

\subsection{$M$-structures related to affine Dynkin diagrams and quiver representations}

We recall that the quiver is just a directed graph $Q=(Ver,E),$
where $Ver$ is a finite set of vertices and $E$ is a finite set of
arrows between them. If $a\in E$ is an arrow, then $t_a$ and $h_a$
denote its tail and its head, respectively. Note that loops and
several arrows with the same tail and head are allowed. A
representation of the quiver $Q$ is a set of vector spaces $L_x$
attached to each vertex $x\in Ver$ and linear maps $f_a: L_{t_a}\to
L_{h_a}$ attached to each arrow $a\in E$. The set of natural numbers
$\mbox{dim} L_x $ attached to each vertex $x\in Ver$ is called the
dimension of the representation. By affine quiver we mean such a
quiver that the corresponding graph is an affine Dynkin diagram of
$ADE$-type.

{\bf Theorem 4.2.} Let $\cal L$ be an $M$-structure with semi-simple
algebras $\cal A$ and $\cal B$ given by (\ref{AB}). Then there
exists a representation of an affine Dynkin quiver such that:

{\bf 1.} There is an one-to-one correspondence between the set of
vector spaces attached to vertices of the quiver and the set of
vector spaces $\{V_1,...,V_r,W_1,...,W_s\}$. Each vector space from
this set is attached to only one vertex.

{\bf 2.} For any $a\in E$ the space attached to its tail $t_a$ is
some of $V_i$ and the space attached to its head $h_a$ is some of
$W_j$.

{\bf 3.} $\cal L$ as $\cal A\otimes\cal B$-module is isomorphic to
$\oplus_{a\in E}V^{\star}_{t_a}\otimes W_{h_a}$.

{\bf 4.} The vector $(\dim V_1,...,\dim V_r,\dim W_1,...,\dim W_s)$
is the minimal imaginary positive root of the Dynkin diagram.

{\bf 5.} The element $1\in{\cal L}=\oplus_{a\in E}Hom(V_{t_a},
W_{h_a})$ is just $\sum_{a\in E}f_a,$ where $f_a$ is the linear map
attached to the arrow $a$.

{\bf Proof.} In the Theorem 4.1 we have already constructed  the
affine Dynkin diagram corresponding to ${\cal L}$ with vector spaces
$\{V_1,...,V_r,W_1,...,W_s\}$ attached to the vertices.  Note that
each edge of this affine Dynkin diagram links some linear spaces
$V_i$ and $W_j$. By definition, the direction of this edge is from
$V_i$ to $W_j$. The decomposition of the element $1\in{\cal
L}=\oplus_{1\le i\le r,1\le j\le s}(V^{\star}_i\otimes
W_j)^{a_{i,j}}$ defines the element from $V^{\star}_i\otimes W_j$.
Since $V^{\star}_i\otimes W_j=Hom(V_i,W_j),$ we obtain a
representation of the quiver. We know already that $J=(\dim
V_1,...,\dim V_r,\dim W_1,...,\dim W_s)$ is an imaginary positive
root. It is easy to see that if it is not minimal, then $\dim{\cal
A\cap\cal B}>1$.

Now we can use known classification of representations of affine
quivers \cite{quivers1,quivers2,quivers3} to describe the
corresponding $M$-structures. Note that each vertex of our quiver
can not be a tail of one arrow and a head of another arrow at the
same time. Given a representation of such a quiver, it remains to
construct an embedding $\cal A\to\cal L$, $\cal B\to\cal L$ and a
scalar product $(\cdot,\cdot)$ on the space $\cal L$. We can
construct the embedding $\cal A\to\cal L$, $\cal B\to\cal L$ by the
formula $a\to 1a$, $b\to b1$ for $a\in\cal A$, $b\in\cal B$ whenever
we know the element $1\in\cal L$. After that it is not difficult to
construct the scalar product.

{\bf Example.} Consider the case $\tilde A_{2 k-1}.$ We have $\dim
V_i=\dim W_i=1$ for $1\le i\le k$. Let $\{v_i\}$ be a basis of
$V^{\star}_i$ and $\{w_i\}$ be a basis of $W_i$. Let $\{e_i\}$ be a
basis of $End(V_i)$ such that $v_ie_i=v_i$ and $\{f_i\}$ be a basis
of $End(W_i)$ such that $f_iw_i=w_i$. A generic element $1\in\cal L$
in a suitable  basis in $V_i$, $W_i$ can be written in the form
$1=\sum_{1\le i\le k}(v_i\otimes w_i+\lambda v_{i+1}\otimes w_i),$
where index $i$ is taken modulo $k$ and $\lambda\in\C$ is a generic
complex number. The embedding $\cal A\to\cal L$, $\cal B\to\cal L$
is the following: $e_i\to 1e_i=v_i\otimes w_i+\lambda v_i\otimes
w_{i-1},\quad f_i\to f_i1=v_i\otimes w_i+\lambda v_{i+1}\otimes
w_i$. It is clear that the vector space ${\cal A\cap\cal B}$ is
spanned by the vector $\sum_i(v_i\otimes w_i+\lambda v_i\otimes
w_{i-1})$ and that the algebra $\cal A\cap\cal B$ is isomorphic to
$\C$.

Let $Q=(Ver,E)$ be an affine quiver and $\rho$ be its representation
constructed by a given $M$-structure $\cal L$ with semi-simple
algebras $\cal A$ and $\cal B$. Let $Ver=Ver_t\sqcup Ver_h,$ where
$Ver_t$ is the set of tails and $Ver_h$ is the set of heads of
arrows. We have $\rho: x\to V_x, y\to W_y, a\to f_a$ for $x\in
Ver_t$, $y\in Ver_h$ and $a\in E$. It turns out that representations
of the algebra $U(\cal L)$ can also be described in terms of
representations of the quiver $Q$.

{\bf Theorem 4.3.} Suppose we have a representation of the algebra
$U(\cal L)$ in a linear space $N;$ then there exists a
representation $\tau: x\to N_x, a\to \phi_a;\,\,\,  x\in Ver, a\in
E$ of the quiver $Q$ such that

{\bf 1.} The restriction of the representation on the subalgebra
${\cal A}\subset U(\cal L)$ is isomorphic to $\oplus_{x\in
Ver_t}V_x\otimes N_x$.

{\bf 2.} The restriction of the representation on the subalgebra
${\cal B}\subset U(\cal L)$ is isomorphic to $\oplus_{x\in
Ver_h}W_x\otimes N_x$.

{\bf 3.} The formula $f=\sum_{a\in E} f_a\otimes \phi_a$ defines an
isomorphism $f:\oplus_{x\in Ver_t}V_x\otimes N_x\to \oplus_{x\in
Ver_h}W_x\otimes N_x.$

{\bf Proof.} It is known that any representation of the algebra
$End(V)$ has the form $V\otimes S,$  where $S$ is a linear space.
The action is given by $f(v\otimes s)=(fv)\otimes s$. Therefore $N$
has the form $N^a=\oplus_{x\in Ver_t}V_x\otimes N_x$ with respect to
the action of ${\cal A}=\oplus_{1\le i\le r}End(V_i)$  and has the
form $N^b=\oplus_{x\in Ver_h}W_x\otimes N_x$ with respect to the
action of ${\cal B}=\oplus_{1\le j\le s}End(W_j)$ for some linear
spaces $N_x$. Both linear spaces $N^a$ and $N^b$ are isomorphic to
$N$. Thus we have linear spaces $N_x$ attached to each $x\in Ver$
and isomorphism $f:\oplus_{x\in Ver_t}V_x\otimes N_x\to \oplus_{x\in
Ver_h}W_x\otimes N_x$. Let $f=\sum_{x,y\in Ver}f_{x,y}$. It is easy
to see that $f_{x,y}=0$ if $x$ and $y$ are not linked by arrow and
$f_{x,y}=f_a\otimes\phi_a$ for some $\phi_a$ if $x=t_a$, $y=h_a$.
Here $f_a$ is defined by Theorem 4.2 (see the property 5). This
gives us a linear map $\phi_a$ attached to each arrow $a\in E$.

{\bf Remark 1.} It is clear that all statements of this section are
valid for weak $M$-structures with semi-simple algebras $\cal A$ and
$\cal B$. However, it is possible to check that any such a weak
$M$-structure has a quadratic central element $K$ and therefore is
an $M$-structure.

{\bf Remark 2.} It is clear from the Theorem 4.3 (see property {\bf
3}) that
$$\dim N=\sum_{x\in Ver_t}m_x\dim N_x=\sum_{x\in Ver_h}n_x\dim N_x.$$
Moreover, if the representation $\tau$ is decomposable then the
representation of $U(\cal L)$ is also decomposable. Therefore,
$\dim\tau$ must be a positive root for indecomposable
representation. If this root is real, then the representation does
not depend on parameters and corresponds to some special value of
$K$. If this root is imaginary, then the representation depends on
one parameter and the action of $K$ also depends on this parameter.
In Appendix we describe these representations for imaginary roots
explicitly.

\section{Appendix}

In this appendix we give explicit formulas for $M$-algebras with
semi-simple algebras $\cal A$ and $\cal B$ based on known
classification results on affine quiver representations. We give
also formulas for the operator $R$ with values in $End(U(\cal L))$.
Note that $K=R(1)$. It turns out that in all cases
\begin{equation}\label{s}
S_{\lambda}=1+\lambda R.
\end{equation}
Moreover, the operator $R$ satisfies some polynomial equation of
degree 3 in the case $\tilde A_{2 k-1}$ and degree 4 in other cases.
Using these equations one can define $(v+R)^{-1}$ with values in the
localization $\C(K)\otimes U(\cal L),$ where $\C(K)$ is the field of
rational functions in $K$. Formula (\ref{ybsol1}) gives us the
corresponding universal $r$-matrix with values in $\C(K)\otimes
U(\cal L)$. For any representation of $U(\cal L)$ in a vector space
$N$ the image of this $r$-matrix is an $r$-matrix with values in
$End(N)$.

{\bf The case $\tilde A_{2 k-1}$}

The algebras $\cal A$ and $\cal B$ have bases $\{e_i; i\in\Z/k\Z\}$
and $\{f_i; i\in\Z/k\Z\}$ correspondingly with multiplications
\begin{equation}\label{mul}
\mathcal{}e_ie_j=\delta_{i,j}e_i, \qquad f_if_j=\delta_{i,j}f_i.
\end{equation}
The $M$-algebra $U(\cal L)$ is generated by
$e_1,...,e_k,f_1,...,f_k$ with defining relations $(\ref{mul})$ and
$$e_1+...+e_k=f_1+...+f_k=1,$$
$$f_ie_j=0, \quad j-i\ne 0,1.$$
The operator $R$ can be written in the form:
$$R(x)=\sum_{1\le i\le j\le k-1}e_ixf_j+f_ke_kx.$$
This operator satisfies the following equation
$$KR(x)-(K+1)R^2(x)+R^3(x)=0.$$
From this equation and (\ref{s}) we obtain
$$(v+R)^{-1}(x)=\frac{1}{v}x+\frac{1}{v(v+1)}(v+K)^{-1}(R^2(x)-(1+v+K)R(x))$$
and $r$-matrix is given by (\ref{ybsol1}).

For each generic value of $K$ the algebra $U(\cal L)$ has the
following irreducible representation $V$. There are two bases
$\{v_i; i\in\Z/k\Z\}$ and $\{w_i; i\in\Z/k\Z\}$ of the space $V$
such that
$$e_iv_j=\delta_{i,j}v_i, \quad f_iw_j=\delta_{i,j}w_i, \quad
v_i=w_i-tw_{i-1}, \quad i,j\in\Z/k\Z.$$ Here $t\in\C$ is a parameter
of representation. In this representation $K$ acts as multiplication
by $1/(1-t^k)$.

{\bf The case $\tilde D_{2 k}$}

The algebra ${\cal A}\cong\C\oplus\C\oplus
(Mat_2)^{k-2}\oplus\C\oplus\C$ has a basis $\{e_1, e_2,e_{2k},
e_{2k+1}, e_{2\alpha,i,j}; 2\le\alpha\le k-1, 1\le i,j\le 2\}$ with
multiplication
\begin{equation}\label{mul3}
e_{\alpha}e_{\beta}=\delta_{\alpha,\beta}e_{\beta}, \quad
e_{\alpha}e_{\beta,i,j}=e_{\beta,i,j}e_{\alpha}=0, \quad
e_{\alpha,i,j}e_{\beta,i^{\prime},j^{\prime}}=\delta_{\alpha,\beta}\delta_{j,i^{\prime}}e_{\alpha,i,j^{\prime}}.
\end{equation}
The algebra ${\cal B}\cong (Mat_2)^{k-1}$ has a basis $\{
e_{2\alpha-1,i,j}; 2\le\alpha\le k, 1\le i,j\le 2\}$ with
multiplication
\begin{equation}\label{mul4}
e_{\alpha,i,j}e_{\beta,i^{\prime},j^{\prime}}=\delta_{\alpha,\beta}\delta_{j,i^{\prime}}e_{\alpha,i,j^{\prime}}.
\end{equation}
The $M$-algebra $U(\cal L)$ is generated by $e_1,e_2,e_{2k},
e_{2k+1}, e_{\alpha,i,j}; 3\le\alpha\le 2k-1, 1\le i,j\le 2$ with
defining relations $(\ref{mul3}), (\ref{mul4})$ and
$$e_1+e_2+e_{2k}+e_{2k+1}+\sum_{2\le\alpha\le k-1, 1\le i\le
2}e_{2\alpha,i,i}=\sum_{2\le\alpha\le k, 1\le i\le
2}e_{2\alpha-1,i,i}=1,$$
$$e_{2\alpha-1,i,j}e_{\beta}=0, \quad 2<\alpha<k, \quad \beta=1,2,2k,2k+1,$$
$$e_{2\alpha-1,i,j}e_{2\beta,i^{\prime},j^{\prime}}=0,  \quad
\alpha\ne\beta, \beta+1,$$
$$e_{3,1,2}e_1=e_{3,2,2}e_1=e_{3,1,1}e_2=e_{3,2,1}e_2=0,$$
$$e_{2\alpha-1,i,j}e_{2\alpha,i^{\prime},j^{\prime}}=e_{2\alpha+1,i,j}e_{2\alpha,i^{\prime},j^{\prime}}=0, \quad
j\ne i^{\prime},$$
$$e_{2\alpha-1,i,1}e_{2\alpha,1,j}=e_{2\alpha-1,i,2}e_{2\alpha,2,j},
\quad
e_{2\alpha+1,i,1}e_{2\alpha,1,j}=e_{2\alpha+1,i,2}e_{2\alpha,2,j},$$
$$e_{2k-1,1,1}e_{2k}=e_{2k-1,1,2}e_{2k}, \quad e_{2k-1,2,1}e_{2k}=e_{2k-1,2,2}e_{2k},$$
$$e_{2k-1,1,2}e_{2k+1}=\lambda e_{2k-1,1,1}e_{2k+1}, \quad e_{2k-1,2,2}e_{2k+1}=\lambda e_{2k-1,2,1}e_{2k+1}.$$
The operator $R$ can be written in the form:
$$R(x)=
\sum_{1\le\alpha\le k-1}(\lambda e_1 x e_{2\alpha+1,2,2}-\lambda e_1
x e_{2\alpha+1,2,1}+e_{2} x e_{2\alpha+1,1,1}- e_2 x
e_{2\alpha+1,1,2}+e_{2k} x e_{2\alpha+1,1,1}+$$
$$\lambda
e_{2k} x e_{2\alpha+1,2,2}+\lambda e_{2k+1} x
e_{2\alpha+1,1,1}+\lambda e_{2k+1} x e_{2\alpha+1,2,2}) +$$
$$
\sum_{2\le\alpha\le k-1,\,\,2\le\beta\le k} (\lambda e_{2\alpha,1,1}
x e_{2\beta-1,2,2}+ e_{2\alpha,2,2} x e_{2\beta-1,1,1})-$$
$$\sum_{2\le\alpha<\beta\le k}(\lambda
e_{2\alpha,1,1} x e_{2\beta-1,2,1}+e_{2\alpha,2,2} x
e_{2\beta-1,1,2})+$$
$$\sum_{2\le\beta\le\alpha\le k-1}(\lambda
e_{2\alpha,2,1} x e_{2\beta-1,2,2}+e_{2\alpha,1,2} x
e_{2\beta-1,1,1}) +(1-\lambda)e_{2k-1,2,2}e_{2k+1} x.$$ This
operator satisfies the following equation
$$R^4(x)-(1+\lambda+K)R^3(x)+(\lambda+K+\lambda K)R^2(x)-\lambda KR(x)=0.$$
From this equation and (\ref{s}) we obtain
$$(v+R)^{-1}(x)=-\frac{1}{v}x+\frac{1}{v(v+1)(v+\lambda)}(v+K)^{-1}\Big(R^3(x)-(1+v+\lambda+K)R^2(x)+$$
$$(v^2+\lambda v+v+ \lambda+(1+v+\lambda)K)R(x)\Big)$$ and $r$-matrix is
given by (\ref{ybsol1}).

For each generic value of $K$ the algebra $U(\cal L)$ has the
following irreducible representation $V$ of dimension $4 k-4$. There
are two bases $\{v_1, v_2, v_{2k}, v_{2k+1}, v_{2\alpha,i,j};
2\le\alpha\le k-1, 1\le i,j\le 2\}$ and $\{ v_{2\alpha-1,i,j};
2\le\alpha\le k, 1\le i,j\le 2\}$ of the space $V$ such that
$$e_{\alpha}v_{\beta}=\delta_{\alpha,\beta}v_{\beta}, \quad
\alpha,\beta=1,2,2k,2k+1,
$$ $$
e_{\alpha}v_{2\beta,i,j}=e_{2\beta,i,j}v_{\alpha}=0,\quad
\alpha=1,2,2k,2k+1 \quad 2\le\beta\le k-1, $$ $$
e_{2\alpha,i,j}v_{2\beta,i^{\prime},j^{\prime}}=\delta_{\alpha,\beta}\delta_{j,i^{\prime}}
v_{2\alpha,i,j^{\prime}}, \quad 2\le\alpha,\beta\le k-1,$$
$$
e_{2\alpha-1,i,j}v_{2\beta-1,i^{\prime},j^{\prime}}=\delta_{\alpha,\beta}\delta_{j,i^{\prime}}
v_{2\alpha-1,i,j^{\prime}}, \quad 2\le\alpha,\beta\le k$$ and
$$v_1=v_{3,1,1}, \quad v_2=v_{3,2,2},$$
$$v_{2\alpha,i,j}=v_{2\alpha+1,i,j}-v_{2\alpha-1,i,j}, \quad
2\le\alpha\le k-1, \quad i,j=1,2, $$
$$v_{2k}=v_{2k-1,1,1}+v_{2k-1,2,1}+v_{2k-1,1,2}+v_{2k-1,2,2},$$
$$v_{2k+1}=v_{2k-1,1,1}+\lambda v_{2k-1,2,1}+tv_{2k-1,1,2}+\lambda tv_{2k-1,2,2}.$$
Here $\lambda\in\C$ is a parameter of the algebra $U(\cal L)$ and
$t\in\C$ is a parameter of representation. In this representation
$K$ acts as multiplication by $\mu=\lambda(t-1)/(t-\lambda)$.

{\bf The case $\tilde D_{2 k-1}$}

The algebra ${\cal A}\cong\C\oplus\C\oplus (Mat_2)^{k-2}$ has a
basis $\{e_1, e_2, e_{2\alpha,i,j}; 2\le\alpha\le k-1, 1\le i,j\le
2\}$ with multiplication
\begin{equation}\label{mul1}
e_{\alpha}e_{\beta}=\delta_{\alpha,\beta}e_{\beta}, \quad
e_{\alpha}e_{\beta,i,j}=e_{\beta,i,j}e_{\alpha}=0, \quad
e_{\alpha,i,j}e_{\beta,i^{\prime},j^{\prime}}=\delta_{\alpha,\beta}\delta_{j,i^{\prime}}e_{\alpha,i,j^{\prime}}.
\end{equation}
The algebra ${\cal B}\cong\C\oplus\C\oplus (Mat_2)^{k-2}$ has a
basis $\{e_{2k-1}, e_{2k}, e_{2\alpha-1,i,j}; 2\le\alpha\le k-1,
1\le i,j\le 2\}$ with multiplication
\begin{equation}\label{mul2}
e_{\alpha}e_{\beta}=\delta_{\alpha,\beta}e_{\beta}, \quad
e_{\alpha}e_{\beta,i,j}=e_{\beta,i,j}e_{\alpha}=0, \quad
e_{\alpha,i,j}e_{\beta,i^{\prime},j^{\prime}}=\delta_{\alpha,\beta}\delta_{j,i^{\prime}}e_{\alpha,i,j^{\prime}}.
\end{equation}
The $M$-algebra $U(\cal L)$ is generated by $e_1,e_2,e_{2k-1},
e_{2k}, e_{\alpha,i,j}; 3\le\alpha\le 2k-2, 1\le i,j\le 2$ with
defining relations $(\ref{mul1}), (\ref{mul2})$ and
$$e_{\alpha}e_{\beta}=0, \quad \alpha=2k-1,2k, \quad \beta=1,2,$$
$$e_1+e_2+\sum_{2\le\alpha\le k-1, 1\le i\le
2}e_{2\alpha,i,i}=e_{2k-1}+e_{2k}+\sum_{2\le\alpha\le k-1, 1\le i\le
2}e_{2\alpha-1,i,i}=1,$$
$$e_{2\alpha-1,i,j}e_{\beta}=0, \quad \alpha>2, \quad \beta=1,2,$$
$$e_{2\alpha-1,i,j}e_{2\beta,i^{\prime},j^{\prime}}=0,  \quad
\alpha\ne\beta, \beta+1,$$
$$e_{\alpha}e_{2\beta,i,j}=0, \quad \beta<k-1, \quad \alpha=2k-1, 2k,$$
$$e_{3,1,2}e_1=e_{3,2,2}e_1=e_{3,1,1}e_2=e_{3,2,1}e_2=0,$$
$$e_{2\alpha-1,i,j}e_{2\alpha,i^{\prime},j^{\prime}}=e_{2\alpha+1,i,j}e_{2\alpha,i^{\prime},j^{\prime}}=0, \quad
j\ne i^{\prime},$$
$$e_{2\alpha-1,i,1}e_{2\alpha,1,j}=e_{2\alpha-1,i,2}e_{2\alpha,2,j},
\quad
e_{2\alpha+1,i,1}e_{2\alpha,1,j}=e_{2\alpha+1,i,2}e_{2\alpha,2,j},$$
$$e_{2k-1}e_{2k-2,1,1}=e_{2k-1}e_{2k-2,2,1}, \quad e_{2k-1}e_{2k-2,1,2}=e_{2k-1}e_{2k-2,2,2}$$
$$e_{2k}e_{2k-2,2,1}=\lambda e_{2k}e_{2k-2,1,1}, \quad e_{2k}e_{2k-2,2,2}=\lambda e_{2k}e_{2k-2,1,2}.$$
The operator $R$ can be written in the form:
$$R(x)=(\lambda-1)e_1xe_{2k-1}+\sum_{2\le\alpha\le
k-1}((\lambda-1)e_1xe_{2\alpha-1,2,2}+(\lambda-1)e_{2\alpha,1,1}xe_{2k-1}-\lambda
e_2xe_{2\alpha-1,1,2}-$$ $$e_1xe_{2\alpha-1,2,1}+\lambda
e_2xe_{2\alpha-1,2,2}+\lambda e_1xe_{2\alpha-1,1,1})+\sum_{2\le
\alpha,\beta\le k-1}((\lambda-1)e_{2\alpha,1,1}xe_{2\beta-1,2,2}+$$
$$\lambda
e_{2\alpha,1,1}xe_{2\beta-1,1,1}+\lambda
e_{2\alpha,2,2}xe_{2\beta-1,2,2})+\sum_{2\le\beta\le\alpha\le
k-1}(\lambda
e_{2\alpha,1,2}xe_{2\beta-1,1,1}+e_{2\alpha,2,1}xe_{2\beta-1,2,2})-$$
$$\sum_{2\le\alpha<\beta\le k-1}(\lambda
e_{2\alpha,2,2}xe_{2\beta-1,1,2}+e_{2\alpha,1,1}xe_{2\beta-1,2,1})+(\lambda-1)xe_{2k}e_{2k-2,2,2}.$$
This operator satisfies the following equation
$$R^4(x)-R^3(x)(2\lambda-1+K)+R^2(x)(\lambda^2-\lambda-K+2\lambda K)-\lambda(\lambda-1)R(x) K=0.$$
From this equation and (\ref{s}) we obtain
$$(v+R)^{-1}(x)=-\frac{1}{v}x+\frac{1}{v(v+\lambda)(v+\lambda-1)}\Big(R^3(x)-R^2(x)(v+2\lambda-1+K)+R(x)(v^2+2\lambda v+
\lambda^2-v-\lambda+$$$$(v-1+2\lambda)K)\Big)(v+K)^{-1}$$ and
$r$-matrix is given by (\ref{ybsol1}).

For each generic value of $K$ the algebra $U(\cal L)$ has the
following irreducible representation $V$ of dimension $4 k-6$. There
are two bases $\{v_1, v_2, v_{2\alpha,i,j}; 2\le\alpha\le k-1, 1\le
i,j\le 2\}$ and $\{v_{2k-1}, v_{2k}, v_{2\alpha-1,i,j};
2\le\alpha\le k-1, 1\le i,j\le 2\}$ of the space $V$ such that
$$e_{\alpha}v_{\beta}=\delta_{\alpha,\beta}v_{\beta}, \quad
\alpha,\beta=1,2,
$$ $$
e_{\alpha}v_{2\beta,i,j}=e_{2\beta,i,j}v_{\alpha}=0,\quad
\alpha=1,2, \quad 2\le\beta\le k-1, $$ $$
e_{2\alpha,i,j}v_{2\beta,i^{\prime},j^{\prime}}=\delta_{\alpha,\beta}\delta_{j,i^{\prime}}v_{2\alpha,i,j^{\prime}},
\quad 2\le\alpha,\beta\le k-1,$$
$$e_{\alpha}v_{\beta}=\delta_{\alpha,\beta}v_{\beta}, \quad
\alpha,\beta=2k-1,2k,
$$ $$
e_{\alpha}v_{2\beta-1,i,j}=e_{2\beta-1,i,j}v_{\alpha}=0,\quad
\alpha=2k-1,2k, \quad 2\le\beta\le k-1, $$ $$
e_{2\alpha-1,i,j}v_{2\beta-1,i^{\prime},j^{\prime}}=\delta_{\alpha,\beta}\delta_{j,i^{\prime}}v_{2\alpha-1,i,j^{\prime}},
\quad 2\le\alpha,\beta\le k-1$$ and
$$v_1=v_{3,1,1}, \quad v_2=v_{3,2,2},$$
$$v_{2\alpha,i,j}=v_{2\alpha+1,i,j}-v_{2\alpha-1,i,j}, \quad
2\le\alpha<k-1, \quad i,j=1,2, $$
$$v_{2k-2,1,1}=v_{2k-1}+v_{2k}-v_{2k-3,1,1}, \quad v_{2k-2,2,2}=v_{2k-1}+\lambda tv_{2k}-v_{2k-3,2,2},$$
$$ v_{2k-2,1,2}=v_{2k-1}+tv_{2k}-v_{2k-3,1,2}, \quad
v_{2k-2,2,1}=v_{2k-1}+\lambda v_{2k}-v_{2k-3,2,1}.$$ Here
$\lambda\in\C$ is a parameter of the algebra $U(\cal L)$ and
$t\in\C$ is a parameter of representation. In this representation
$K$ acts as multiplication by
$\mu=t\lambda(1-\lambda)/(1-t\lambda)$.

{\bf The case $\tilde E_6$}

The algebra ${\cal A}\cong\C\oplus\C\oplus\C\oplus Mat_3$ has a
basis $\{e_{\alpha}, e_{\beta,\gamma};
\alpha,\beta,\gamma\in\Z/3\Z\}$ with multiplication
\begin{equation}\label{mul5}
e_{\alpha}e_{\beta}=\delta_{\alpha,\beta}e_{\beta}, \quad
e_{\alpha}e_{\beta,\gamma}=e_{\beta,\gamma}e_{\alpha}=0, \quad
e_{\beta,\gamma}e_{\beta^{\prime},\gamma^{\prime}}=\delta_{\gamma,\beta^{\prime}}e_{\beta,\gamma^{\prime}}.
\end{equation}
The algebra ${\cal B}\cong Mat_2\oplus Mat_2\oplus Mat_2$ has a
basis $\{ f_{\alpha,i,j}; \alpha\in\Z/3\Z, 1\le i,j\le 2\}$ with
multiplication
\begin{equation}\label{mul6}
f_{\alpha,i,j}f_{\beta,i^{\prime},j^{\prime}}=\delta_{\alpha,\beta}\delta_{j,i^{\prime}}f_{\alpha,i,j^{\prime}}.
\end{equation}
The $M$-algebra $U(\cal L)$ is generated by $e_{\alpha},
e_{\beta,\gamma}, f_{\alpha,i,j}; \alpha,\beta,\gamma\in\Z/3\Z, 1\le
i,j\le 2$ with defining relations $(\ref{mul5}), (\ref{mul6})$ and
$$\sum_{\alpha\in\Z/3\Z}e_{\alpha}+\sum_{\beta\in\Z/3\Z}e_{\beta,\beta}=
\sum_{1\le\alpha\le 3, 1\le i\le 2}f_{\alpha,i,i}=1,$$
$$f_{\alpha,i,j}e_{\beta}=0, \quad \alpha\ne\beta, \quad f_{\alpha,i,2}e_{\alpha}=\lambda f_{\alpha,i,1}e_{\alpha},$$
$$f_{\alpha,i,j}e_{\alpha+1,\beta}=0, \quad f_{\alpha,i,1}e_{\alpha,\beta}=0, \quad f_{\alpha+1,i,2}e_{\alpha,\beta}=0,$$
$$f_{\alpha,i,2}e_{\alpha,\beta}=f_{\alpha,i,1}e_{\alpha-1,\beta}.$$
The operator $R$ can be written in the form:
$$R(x)=(\lambda^3-1) \Big(e_{1} x f_{1,2,2}+e_{1} x f_{2,1,1}+e_{1} x f_{2,2,2}+e_{1} x f_{3,1,1}
+e_{2} x f_{2,2,2}+e_{2} x f_{3,1,1}+ $$
$$ e_{1,1} x f_{2,2,2}+ e_{1,1} x f_{3,1,1}+e_{3,3} x f_{1,2,2}+e_{3,3} x f_{2,1,1}
+e_{3,3} x f_{2,2,2}+e_{3,3} x f_{3,1,1}\Big)+ $$
$$
\lambda^3 \Big(e_{1,1} x f_{1,2,2}+e_{1,1} x f_{2,1,1}+e_{2,2} x
f_{2,2,2}+e_{2,2} x f_{3,1,1} +e_{3,3} x f_{1,1,1}+e_{3,3} x
f_{3,2,2}\Big)+
$$
$$
\lambda^2 \Big(e_{1,1} x f_{1,1,2}+e_{2,2} x f_{2,1,2}+e_{3,3} x
f_{3,1,2}- e_{1,2} x f_{1,2,2} - e_{1,3} x f_{2,2,1}-
$$
$$ e_{2,1} x f_{3,2,1}-e_{2,3} x f_{2,2,2}- e_{3,1} x f_{3,2,2}-e_{3,2} x f_{1,2,1}\Big)+
$$
$$
\lambda \Big(e_{1,1} x f_{2,2,1}+e_{2,2} x f_{3,2,1}+e_{3,3} x
f_{1,2,1}- e_{1,2} x f_{1,1,2} - e_{1,3} x f_{2,1,1}-
$$
$$ e_{2,1} x f_{3,1,1}-e_{2,3} x f_{2,1,2}- e_{3,1} x
f_{3,1,2}-e_{3,2} x  f_{1,1,1}\Big)+(\lambda^3-1) f_{3,2,2} e_{3} x.
$$
This operator satisfies the following equation
$$R^4(x)-(2\lambda^3-1+K)R^3(x)+(\lambda^6-\lambda^3-K+2\lambda^3 K)R^2(x)-\lambda^3(\lambda^3-1)KR(x) =0.$$
From this equation and (\ref{s}) we obtain
$$(v+R)^{-1}(x)=-\frac{1}{v}x+\frac{1}{v(v+\lambda^3)(v+\lambda^3-1)}(v+K)^{-1}\Big(R^3(x)-(v+2\lambda^3-1+K)R^2(x)+$$$$
(v^2+2\lambda^3 v+
\lambda^6-v-\lambda^3+(v-1+2\lambda^3)K)R(x)\Big)$$ and $r$-matrix
is given by (\ref{ybsol1}).

For each generic value of $K$ the algebra $U(\cal L)$ has the
following irreducible representation $V$. There are two bases
$\{v_{\alpha}, v_{\beta,\gamma}; \alpha,\beta,\gamma\in\Z/3\Z\}$ and
$\{ w_{\alpha,i,j}; \alpha\in\Z/3\Z, 1\le i,j\le 2\}$ of the space
$V$ such that
$$e_{\alpha}v_{\beta}=\delta_{\alpha,\beta}v_{\beta}, \quad
e_{\alpha}v_{\beta,\gamma}=e_{\beta,\gamma}v_{\alpha}=0, \quad
e_{\beta,\gamma}v_{\beta^{\prime},\gamma^{\prime}}=\delta_{\gamma,\beta^{\prime}}v_{\beta,\gamma^{\prime}},$$
$$f_{\alpha,i,j}w_{\beta,i^{\prime},j^{\prime}}=\delta_{\alpha,\beta}\delta_{j,i^{\prime}}w_{\alpha,i,j^{\prime}}$$
and
$$v_{\alpha}=w_{\alpha,1,1}+\lambda w_{\alpha,2,1}+tw_{\alpha,1,2}+\lambda tw_{\alpha,2,2},$$
$$v_{\alpha,\alpha}=w_{\alpha,2,2}, \quad
v_{\alpha+1,\alpha}=w_{\alpha+2,1,1}, \quad
v_{\alpha+2,\alpha}=w_{\alpha+2,2,1}+w_{\alpha,1,2}.$$ Here
$\lambda\in\C$ is a parameter of the algebra $U(\cal L)$ and
$t\in\C$ is a parameter of representation. In this representation
$K$ acts as multiplication by $\displaystyle \mu=\frac{\lambda^3
(\lambda^3-1)}{(\lambda^3-t^3)}$.

{\bf The case $\tilde E_7$}

The algebra ${\cal A}\cong\C\oplus Mat_3\oplus Mat_2\oplus
Mat_3\oplus\C$ has a basis $\{e_1,
e_{2,i_1,j_1},e_{3,i_2,j_2},e_{4,i_3,j_3}, e_5 ; 1\le
i_1,j_1,i_3,j_3\le 3, 1\le i_2,j_2\le 2\}$ with multiplication
\begin{equation}\label{mul7}
e_{\alpha}e_{\beta}=\delta_{\alpha,\beta}e_{\beta}, \quad
e_{\alpha}e_{\beta,i,j}=e_{\beta,i,j}e_{\alpha}=0, \quad
e_{\alpha,i,j}e_{\beta,i^{\prime},j^{\prime}}=\delta_{\alpha,\beta}\delta_{j,i^{\prime}}e_{\alpha,i,j^{\prime}}.
\end{equation}
The algebra ${\cal B}\cong Mat_2\oplus Mat_4\oplus Mat_2$ has a
basis $\{ f_{1,i_1,j_1},f_{2,i_2,j_2},f_{3,i_3,j_3}; 1\le
i_1,j_1,i_3,j_3\le 2, 1\le i_2,j_2\le 4\}$ with multiplication
\begin{equation}\label{mul8}
f_{\alpha,i,j}f_{\beta,i^{\prime},j^{\prime}}=\delta_{\alpha,\beta}\delta_{j,i^{\prime}}f_{\alpha,i,j^{\prime}}.
\end{equation}
The $M$-algebra $U(\cal L)$ is generated by $e_1,
e_{2,i,j},e_{3,i,j},e_{4,i,j}, e_5, f_{1,i,j},f_{2,i,j},f_{3,i,j}$
with defining relations $(\ref{mul7}), (\ref{mul8})$ and
$$e_1+
\sum_{1\le i\le 3}e_{2,i,i}+\sum_{1\le i\le 2}e_{3,i,i}+\sum_{1\le
i\le 3}e_{4,i,i}+e_5= \sum_{1\le i\le 2}f_{1,i,i}+\sum_{1\le i\le
4}f_{2,i,i}+\sum_{1\le i\le 2}f_{3,i,i}=1,$$
$$f_{1,i,j}e_{3,i^{\prime},j^{\prime}}=f_{1,i,j}e_{4,i^{\prime},j^{\prime}}=f_{1,i,j}e_5=f_{2,i,j}e_1=
f_{2,i,j}e_5=f_{3,i,j}e_1=f_{3,i,j}e_{2,i^{\prime},j^{\prime}}=f_{3,i,j}e_{3,i^{\prime},j^{\prime}}=0$$
$$f_{1,i,2}e_1=f_{3,i,2}e_5=0, \quad
f_{1,i,j}e_{2,i^{\prime},j^{\prime}}=f_{3,i,j}e_{4,i^{\prime},j^{\prime}}=0,
\quad j\ne i^{\prime},$$
$$f_{1,i,1}e_{2,1,j}=f_{1,i,2}e_{2,2,j}, \quad
f_{3,i,1}e_{4,1,j}=f_{3,i,2}e_{4,2,j},$$
$$f_{2,i,j}e_{2,i^{\prime},j^{\prime}}=0, \quad j\ne i^{\prime},$$
$$f_{2,i,1}e_{2,1,j}=f_{2,i,2}e_{2,2,j}=f_{2,i,3}e_{2,3,j},$$
$$f_{2,i,j}e_{4,i^{\prime},j^{\prime}}=0, \quad (j,i^{\prime})\notin \{(1,2),(2,1),(4,3)\},$$
$$f_{2,i,1}e_{4,2,j}=f_{2,i,2}e_{4,1,j}=f_{2,i,4}e_{4,3,j},$$
$$f_{2,i,1}e_{3,2,j^{\prime}}=f_{2,i,2}e_{3,1,j}=0,$$
$$f_{2,i,1}e_{3,1,j}=f_{2,i,2}e_{3,2,j}=f_{2,i,3}e_{3,1,j}=f_{2,i,3}e_{3,2,j}=f_{2,i,4}e_{3,1,j}=
\lambda^{-1}f_{2,i,4}e_{3,2,j}.$$ The operator $R$ can be written in
the form:
$$R(x)=\lambda (\lambda-1) (e_{2, 1, 2}xf_{2, 4, 1}+
 e_{2, 3, 2}xf_{2, 4, 3}+ e_{2, 2, 2}xf_{2, 4, 2}-e_{4, 3, 1}xf_{2,
4, 4})+ (\lambda+1) (e_5xf_{3, 2, 2}+ e_{3, 2, 2}xf_{1, 1, 1}+$$
$$ e_5xf_{1, 1, 1}+ e_{2, 2, 2}xf_{1, 1, 1}+ e_5xf_{2, 1, 1}+ e_{3, 2, 2}xf_{2, 1, 1}+
 e_{4, 1, 1}xf_{2, 1, 1}+
e_{4, 1, 1}xf_{1, 1, 1}+ e_{2, 2, 2}xf_{2, 1, 1}+ e_{3, 2, 2}xf_{3,
2, 2}+$$ $$ e_{4, 1, 1}xf_{3, 2, 2}+ e_{2, 2, 2}xf_{3, 2,
2})+(\lambda-1)( e_{4, 1, 2}xf_{2, 3, 2}+ e_{4, 2, 2}xf_{2, 3, 1}+
e_{4, 3, 2}xf_{2, 3, 4}- e_{2, 3, 1}xf_{2, 3, 3})+ \lambda (e_{2, 2,
2}xf_{1, 1, 2}-$$ $$e_{2, 2, 3}xf_{1, 1, 2}+e_{2, 1, 1}xf_{1, 1,
1}+e_{4, 1, 1}xf_{2, 1, 2}+e_{4, 1, 2}xf_{2, 1, 2}+e_5xf_{2, 1, 2}+
e_{4, 3, 1}xf_{2, 1, 4}+e_{4, 3, 2}xf_{2, 1, 4}+e_{4, 1, 1}xf_{1, 1,
2}+$$ $$ e_{4, 1, 2}xf_{1, 1, 2}+e_{2, 2, 2}xf_{3, 1, 1}+e_{2, 2,
2}xf_{2, 2, 2}- e_{4, 3, 1}xf_{2, 2, 4}+e_5xf_{1, 1, 2}+e_{4, 1,
1}xf_{3, 1, 1}- e_{2, 1, 3}xf_{1, 1, 1}-e_{3, 1, 2}xf_{1, 1, 1}+$$
$$e_{4, 3, 3}xf_{2, 1, 1}+e_{4, 2, 2}xf_{1, 1, 1}+ e_{4, 3, 3}xf_{1,
1, 1}+e_{2, 3, 2}xf_{2, 2, 3}- e_{2, 1, 2}xf_{2, 1, 1}+e_5xf_{3, 2,
1}-e_{2, 3, 1}xf_{2, 1, 3}-e_{3, 1, 2}xf_{2, 1, 1}-$$
$$e_{2, 3, 2}xf_{2, 1, 3}+ e_{2, 3, 3}xf_{2, 1, 3}+e_{2, 1, 2}xf_{3,
1, 2}+e_{4, 2, 2}xf_{2, 1, 1}+ e_{2, 2, 2}xf_{2, 3, 3}+e_{2, 3,
3}xf_{2, 3, 3}+e_{2, 1, 2}xf_{2, 2, 1}+ e_{3, 2, 2}xf_{2, 3, 3}+$$
$$e_{4, 1, 1}xf_{2, 3, 3}-e_{2, 1, 2}xf_{3, 2, 2}+e_{4, 3, 3}xf_{2,
4, 4}+e_{2, 2, 1}xf_{1, 1, 2}+e_{4, 2, 1}xf_{3, 1, 2}-e_{3, 1,
2}xf_{3, 2, 2}-e_{4, 2, 1}xf_{3, 2, 2}+ e_{4, 3, 3}xf_{2, 3, 3}+$$
$$e_5xf_{2, 3, 3}+e_{4, 3, 3}xf_{3, 2, 2})+ e_{2, 2,
1}xf_{2, 1, 1}-e_{2, 1, 1}xf_{1, 2, 1}-e_{2, 2, 2}xf_{1, 2, 1}+e_{2,
3, 3}xf_{2, 1, 1}+ e_{3, 2, 1}xf_{2, 1, 1}+e_{4, 1, 2}xf_{2, 1,
1}-$$ $$e_{2, 2, 2}xf_{3, 1, 2}-e_{2, 2, 2}xf_{2, 2, 1}- e_{2, 3,
3}xf_{3, 1, 2}-e_{3, 1, 1}xf_{3, 1, 2}-e_{3, 2, 2}xf_{3, 1, 2}-e_{4,
1, 1}xf_{3, 1, 2}+ e_{2, 2, 1}xf_{2, 3, 3}+e_{3, 2, 1}xf_{2, 3,
3}-$$ $$e_{2, 3, 3}xf_{1, 2, 1}-e_{3, 1, 1}xf_{1, 2, 1}- e_{3, 2,
2}xf_{1, 2, 1}-e_{4, 1, 1}xf_{1, 2, 1}-e_{2, 3, 3}xf_{2, 2, 1}-e_{3,
1, 1}xf_{2, 2, 1}-e_{3, 2, 2}xf_{2, 2, 1}- e_{4, 1, 1}xf_{2, 2,
1}-$$ $$e_{4, 2, 3}xf_{3, 1, 2}-e_{4, 3, 3}xf_{3, 1, 2}-e_5xf_{3, 1,
2}+ e_{2, 2, 1}xf_{3, 2, 2}+e_5xf_{2, 4, 4}+e_{2, 2, 1}xf_{3, 1, 1}-
e_{4, 2, 2}xf_{2, 2, 1}-e_{4, 3, 3}xf_{2, 2, 1}-$$
$$e_5xf_{2, 2, 1}-e_{4, 2, 2}xf_{1, 2, 1}- e_{4, 3, 3}xf_{1,
2, 1}-e_5xf_{1, 2, 1}+e_{4, 1, 2}xf_{2, 3, 3}+e_{2, 3, 3}xf_{3, 2,
2}+ e_{3, 2, 1}xf_{3, 2, 2}+e_{4, 1, 2}xf_{3, 2, 2}+$$
$$e_{2, 3, 3}xf_{2, 4, 4}+e_{3, 2, 1}xf_{2, 4, 4}+ e_{3, 2, 2}xf_{2,
4, 4}+e_{4, 1, 1}xf_{2, 4, 4}+e_{4, 1, 2}xf_{2, 4, 4}+e_{4, 1,
2}xf_{3, 1, 1}-e_{4, 1, 3}xf_{3, 1, 1}+ e_{2, 2, 1}xf_{2, 2, 2}+$$
$$e_{2, 2, 1}xf_{1, 1, 1}+e_{2, 3, 3}xf_{1, 1, 1}-e_{4, 3, 2}xf_{2,
2, 4}+ e_{4, 3, 3}xf_{2, 2, 4}+e_{2, 2, 1}xf_{2, 4, 4}+e_{2, 2,
2}xf_{2, 4, 4}+e_{2, 3, 1}xf_{2, 2, 3}+ e_{3, 2, 1}xf_{1, 1, 1}+$$
$$e_{4, 1, 2}xf_{1, 1, 1}+f_{1, 2, 1} e_1x.
$$

For each generic value of $K$ the algebra $U(\cal L)$ has the
following irreducible representation $V$. There are two bases
$\{v_1, v_{2,i_1,j_1},v_{3,i_2,j_2},v_{4,i_3,j_3}, v_5 ; 1\le
i_1,j_1,i_3,j_3\le 3, 1\le i_2,j_2\le 2\}$ and $\{
w_{1,i_1,j_1},w_{2,i_2,j_2},w_{3,i_3,j_3}; 1\le i_1,j_1,i_3,j_3\le
2, 1\le i_2,j_2\le 4\}$ of the space $V$ such that
$$e_{\alpha}v_{\beta}=\delta_{\alpha,\beta}v_{\beta}, \quad
e_{\alpha}v_{\beta,i,j}=e_{\beta,i,j}v_{\alpha}=0, \quad
e_{\alpha,i,j}v_{\beta,i^{\prime},j^{\prime}}=\delta_{\alpha,\beta}\delta_{j,i^{\prime}}v_{\alpha,i,j^{\prime}},$$
$$f_{\alpha,i,j}w_{\beta,i^{\prime},j^{\prime}}=\delta_{\alpha,\beta}\delta_{j,i^{\prime}}w_{\alpha,i,j^{\prime}}$$
and
$$v_1=w_{1,1,1}, \quad v_5=w_{3,1,1},$$
$$v_{2,i,j}=w_{1,i,j}+w_{2,i,j}, \quad
v_{2,i,3}=w_{2,i,3}$$ $$ v_{2,3,j}=w_{2,3,j}, \quad
v_{2,3,3}=w_{2,3,3}, \quad i,j\le 2,$$
$$v_{4,i,j}=w_{3,i,j}+w_{2,\sigma_i,\sigma_j}, \quad
v_{4,i,3}=w_{2,\sigma_i,4}, \quad v_{4,3,j}=w_{2,4,\sigma_j},$$
$$ v_{4,3,3}=w_{2,4,4}, \quad i,j\le 2, \quad \sigma_1=2, \quad \sigma_2=1,$$
$$v_{3,1,1}=w_{2,1,1}+w_{2,1,3}+w_{2,3,1}+w_{2,3,3}+w_{2,1,4}+w_{2,3,4}+w_{2,4,1}+w_{2,4,3}+w_{2,4,4},$$
$$v_{3,1,2}=w_{2,1,2}+w_{2,1,3}+w_{2,3,2}+w_{2,3,3}+w_{2,4,2}+w_{2,4,3}+t(w_{2,1,4}+w_{2,3,4}+w_{2,4,4}),$$
$$v_{3,2,1}=w_{2,2,1}+w_{2,2,3}+w_{2,2,4}+w_{2,3,1}+w_{2,3,3}+w_{2,3,4}+\lambda(w_{2,4,1}+w_{2,4,3}+w_{2,4,4}),$$
$$v_{3,2,2}=w_{2,2,2}+w_{2,2,3}+w_{2,3,2}+w_{2,3,3}+\lambda(w_{2,4,2}+w_{2,4,3})+t(w_{2,2,4}+w_{2,3,4})+\lambda
tw_{2,4,4}.$$ Here $\lambda\in\C$ is a parameter of the algebra
$U(\cal L)$ and $t\in\C$ is a parameter of representation. In this
representation $K$ acts as $\displaystyle
\frac{\lambda(t-1)}{t-\lambda}$.

{\bf The case $\tilde E_8$}

The algebra ${\cal A}\cong Mat_2\oplus Mat_6\oplus Mat_4\oplus
Mat_2$ has a basis $\{e_{1,i_1,j_1},e_{2,i_2,j_2},e_{3,i_3,j_3},
e_{4,i_4,j_4} ; 1\le i_1,j_1,i_4,j_4\le 2, 1\le i_2,j_2\le 6, 1\le
i_3,j_3\le 4\}$ with multiplication
\begin{equation}\label{mul9}
e_{\alpha,i,j}e_{\beta,i^{\prime},j^{\prime}}=\delta_{\alpha,\beta}\delta_{j,i^{\prime}}e_{\alpha,i,j^{\prime}}.
\end{equation}
The algebra ${\cal B}\cong Mat_4\oplus Mat_3\oplus Mat_5\oplus
Mat_3\oplus\C$ has a basis $\{
f_{1,i_1,j_1},f_{2,i_2,j_2},f_{3,i_3,j_3},f_{4,i_4,j_4},f_5; 1\le
i_1,j_1\le 4, 1\le i_2,j_2,i_4,j_4\le 3, 1\le i_3,j_3\le 5\}$ with
multiplication
\begin{equation}\label{mul10}
f_5^2=f_5, \quad f_5f_{\alpha,i,j}=f_{\alpha,i,j}f_5=0, \quad
f_{\alpha,i,j}f_{\beta,i^{\prime},j^{\prime}}=\delta_{\alpha,\beta}\delta_{j,i^{\prime}}f_{\alpha,i,j^{\prime}}.
\end{equation}
The $M$-algebra $U(\cal L)$ is generated by
$e_{1,i,j},e_{2,i,j},e_{3,i,j},
e_{4,i,j},f_{1,i,j},f_{2,i,j},f_{3,i,j},f_{4,i,j},f_5$ with defining
relations $(\ref{mul9}), (\ref{mul10})$ and
$$\sum_{1\le i\le 2}e_{1,i,i}+\sum_{1\le i\le 6}e_{2,i,i}+\sum_{1\le
i\le 4}e_{3,i,i}+\sum_{1\le i\le 2}e_{4,i,i}= \sum_{1\le i\le
4}f_{1,i,i}+\sum_{1\le i\le 3}f_{2,i,i}+\sum_{1\le i\le
5}f_{3,i,i}+\sum_{1\le i\le 3}f_{4,i,i}+f_5=1,$$
$$f_{1,i,j}e_{3,i^{\prime},j^{\prime}}=f_{1,i,j}e_{4,i^{\prime},j^{\prime}}=f_{2,i,j}e_{1,i^{\prime},j^{\prime}}=
f_{2,i,j}e_{3,i^{\prime},j^{\prime}}=f_{2,i,j}e_{4,i^{\prime},j^{\prime}}=f_{3,i,j}e_{1,i^{\prime},j^{\prime}}=
f_{3,i,j}e_{4,i^{\prime},j^{\prime}}=$$
$$f_{4,i,j}e_{1,i^{\prime},j^{\prime}}=f_{4,i,j}e_{2,i^{\prime},j^{\prime}}=
f_5e_{1,i,j}=f_5e_{2,i,j}=f_5e_{3,i,j}=0,$$
$$f_5e_{4,2,i}=0, \quad f_{4,i,1}e_{4,1,j}=f_{4,i,2}e_{4,2,j}, \quad f_{4,i,j}e_{4,i^{\prime},j^{\prime}}=0, \quad j\ne
i^{\prime},$$ $$
f_{4,i,1}e_{3,1,j}=f_{4,i,2}e_{3,2,j}=f_{4,i,3}e_{3,3,j}, \quad
f_{4,i,j}e_{3,i^{\prime},j^{\prime}}=0, \quad j\ne i^{\prime},$$
$$f_{3,i,1}e_{3,1,j}=f_{3,i,2}e_{3,2,j}=f_{3,i,3}e_{3,3,j}=f_{3,i,4}e_{3,4,j}, \quad
f_{3,i,j}e_{3,i^{\prime},j^{\prime}}=0, \quad j\ne i^{\prime},$$
$$f_{3,i,1}e_{2,1,j}=f_{3,i,2}e_{2,2,j}=f_{3,i,3}e_{2,3,j}=f_{3,i,4}e_{2,4,j}=f_{3,i,5}e_{2,5,j}, \quad
f_{3,i,j}e_{2,i^{\prime},j^{\prime}}=0, \quad j\ne i^{\prime},$$
$$f_{1,i,1}e_{1,1,j}=f_{1,i,2}e_{1,2,j}, \quad
f_{1,i,j}e_{1,i^{\prime},j^{\prime}}=0, \quad j\ne i^{\prime},$$
$$f_{1,i,1}e_{2,1,j}=f_{1,i,2}e_{2,2,j}=f_{1,i,3}e_{2,4,j}=f_{1,i,4}e_{2,6,j}, \quad
f_{1,i,j}e_{2,i^{\prime},j^{\prime}}=0, \quad (j,i^{\prime})\notin
\{(1,1),(2,2),(3,4),(4,6)\},$$
$$f_{2,i,1}e_{2,2,j}=f_{2,i,1}e_{2,4,j}=f_{2,i,1}e_{2,5,j}=f_{2,i,2}e_{2,1,j}=f_{2,i,2}e_{2,3,j}=f_{2,i,2}e_{2,4,j}=
f_{2,i,3}e_{2,1,j}=f_{2,i,3}e_{2,6,j}=$$
$$\lambda^{-1}f_{2,i,3}e_{2,5,j},
\quad f_{2,i,j}e_{2,i^{\prime},j^{\prime}}=0, \quad
(j,i^{\prime})\notin
\{(1,2),(1,4),(1,5),(2,1),(2,3),(2,4),(3,1),(3,6),(3,5)\}.$$ The
operator $R$ can be written in the form:
$$R(x)=(\lambda+1) (e_{2, 1, 1} xf_{1, 3, 3}+e_{3, 1, 1}
xf_{1, 3, 3}+e_{2, 6, 6} xf_{1, 3, 3}+e_{2, 3, 3} xf_{2, 1, 1}+e_{2,
1, 1} xf_{3, 4, 4}+ e_{4, 1, 1} xf_{1, 3, 3}+$$ $$e_{2, 1, 1} xf_{2,
1, 1}+e_{1, 1, 1} xf_{2, 1, 1}+e_{3, 1, 1} xf_{2, 1, 1}+e_{2, 6, 6}
xf_{2, 1, 1}+e_{1, 1, 1} xf_{3, 4, 4}+e_{4, 1, 1} xf_{3, 4, 4}+
e_{3, 3, 3} xf_{2, 1, 1}+e_{4, 1, 1} xf_{2, 1, 1}+$$ $$e_{2, 6, 6}
xf_{3, 4, 4}+e_{3, 1, 1} xf_{3, 5, 5}+e_{3, 1, 1} xf_{3, 4, 4}+e_{2,
3, 3} xf_{3, 5, 5}+e_{1, 1, 1} xf_{3, 5, 5}+ e_{2, 1, 1} xf_{3, 5,
5}+e_{3, 3, 3} xf_{3, 5, 5}+e_{4, 1, 1} xf_{3, 5, 5}+$$ $$e_{2, 6,
6} xf_{1, 2, 2}+e_{3, 3, 3} xf_{1, 2, 2}+e_{1, 1, 1} xf_{1, 2,
2}+e_{2, 1, 1} xf_{1, 2, 2}+ e_{3, 1, 1} xf_{1, 2, 2}+e_{2, 3, 3}
xf_{1, 2, 2}+e_{4, 1, 1} xf_{1, 2, 2}+e_{2, 1, 1} xf_{4, 2, 2}+$$
$$e_{1, 1, 1} xf_{4, 2, 2}+e_{1, 1, 1} xf_{3, 2, 2}+e_{2, 1, 1}
xf_{3, 2, 2}+ e_{2, 6, 6} xf_{4, 2, 2}+e_{3, 1, 1} xf_{4, 2,
2}+e_{2, 3, 3} xf_{4, 2, 2}+e_{1, 1, 1} xf_{1, 3, 3}+e_{3, 1, 1}
xf_{3, 2, 2}+$$ $$e_{2, 6, 6} xf_{3, 2, 2}+e_{2, 3, 3} xf_{3, 2, 2}+
e_{3, 3, 3} xf_{4, 2, 2}+e_{4, 1, 1} xf_{4, 2, 2}+e_{3, 3, 3} xf_{3,
2, 2}+e_{4, 1, 1} xf_{3, 2, 2})+$$ $$(\lambda-1) (e_{2, 3, 1} xf_{4,
3, 3}-e_{2, 4, 4} xf_{1, 2, 3}-e_{3, 4, 4} xf_{1, 2, 3}+ e_{3, 3, 4}
xf_{4, 2, 3}+e_{3, 3, 1} xf_{4, 3, 3}+e_{2, 1, 4} xf_{3, 2, 1}+e_{2,
1, 4} xf_{4, 2, 1}-$$ $$e_{3, 4, 4} xf_{3, 2, 4}-e_{2, 6, 4} xf_{1,
2, 4}-e_{3, 4, 1} xf_{3, 3, 4}-e_{3, 2, 1} xf_{3, 3, 2}+ e_{2, 3, 1}
xf_{3, 3, 3}+e_{2, 5, 4} xf_{3, 2, 5}+e_{2, 5, 1} xf_{3, 3, 5}-e_{4,
1, 1} xf_{3, 3, 1}-$$ $$e_{4, 2, 1} xf_{3, 3, 2}+e_{3, 1, 4} xf_{4,
2, 1}-e_{3, 1, 1} xf_{3, 3, 1}+e_{2, 3, 4} xf_{4, 2, 3}+ e_{3, 2, 4}
xf_{4, 2, 2}+e_{2, 2, 4} xf_{3, 2, 2}-e_{4, 1, 1} xf_{4, 3, 1}+e_{2,
2, 4} xf_{4, 2, 2}+$$ $$e_{2, 3, 4} xf_{3, 2, 3}-e_{4, 2, 1} xf_{4,
3, 2})+\lambda (e_{1, 1, 1} xf_{1, 4, 1}+e_{3, 4, 3} xf_{1, 1,
3}-e_{3, 4, 2} xf_{1, 1, 3}+e_{3, 4, 1} xf_{1, 1, 3}- e_{2, 4, 2}
xf_{1, 2, 3}+ e_{2, 1, 1} xf_{1, 1, 1}+$$ $$ e_{1, 1, 1} xf_{1, 1,
2}+ e_{3, 1, 3} xf_{1, 1, 1}- e_{2, 4, 5} xf_{1, 2, 3}- e_{2, 1, 2}
xf_{1, 1, 1}- e_{3, 4, 2} xf_{1, 2, 3}+ e_{2, 4, 2} xf_{1, 3, 3}-
e_{2, 1, 2} xf_{1, 3, 3}+ e_{2, 1, 3} xf_{1, 1, 1}-$$ $$ e_{1, 1, 2}
xf_{1, 3, 3}- e_{2, 6, 4} xf_{1, 3, 3}+ e_{2, 4, 4} xf_{1, 3, 3}+
e_{3, 3, 2} xf_{4, 2, 3}+ e_{2, 3, 3} xf_{1, 3, 3}+ e_{2, 3, 5}
xf_{3, 1, 3}+ e_{2, 5, 5} xf_{1, 3, 3}+ e_{3, 4, 4} xf_{1, 3, 3}+$$
$$ e_{3, 3, 3} xf_{1, 3, 3}+ e_{3, 4, 2} xf_{1, 3, 3}- e_{3, 1, 2}
xf_{1, 3, 3}+ e_{2, 1, 1} xf_{1, 4, 3}+ e_{1, 2, 1} xf_{1, 4, 3}+
e_{1, 1, 1} xf_{1, 4, 3}- e_{4, 1, 2} xf_{1, 3, 3}+ e_{2, 4, 2}
xf_{1, 4, 3}-$$ $$ e_{2, 4, 1} xf_{1, 4, 3}+ e_{2, 2, 1} xf_{1, 4,
3}- e_{1, 1, 2} xf_{4, 3, 3}+ e_{3, 1, 1} xf_{1, 4, 3}+ e_{2, 6, 6}
xf_{1, 4, 3}- e_{2, 1, 2} xf_{4, 3, 3}+ e_{2, 5, 5} xf_{1, 4, 3}+
e_{2, 3, 5} xf_{1, 4, 3}+$$ $$ e_{2, 4, 4} xf_{1, 4, 3}+ e_{2, 3, 3}
xf_{4, 3, 3}+ e_{3, 2, 1} xf_{1, 4, 3}- e_{2, 1, 5} xf_{1, 1, 1}+
e_{3, 4, 2} xf_{1, 4, 3}- e_{3, 4, 1} xf_{1, 4, 3}+ e_{3, 4, 1}
xf_{3, 1, 4}+ e_{3, 1, 3} xf_{3, 1, 1}+$$ $$ e_{4, 2, 1} xf_{1, 4,
3}- e_{2, 3, 5} xf_{4, 3, 3}+ e_{4, 1, 1} xf_{1, 4, 3}+ e_{3, 4, 4}
xf_{1, 4, 3}- e_{2, 6, 4} xf_{4, 3, 3}- e_{2, 6, 2} xf_{1, 1, 4}+
e_{3, 4, 3} xf_{3, 1, 4}- e_{3, 4, 2} xf_{3, 1, 4}+$$ $$ e_{3, 1, 1}
xf_{1, 1, 1}+ e_{2, 6, 1} xf_{1, 1, 4}+ e_{3, 1, 1} xf_{3, 3, 5}+
e_{2, 6, 3} xf_{1, 1, 4}- e_{3, 1, 2} xf_{4, 3, 3}- e_{3, 1, 2}
xf_{1, 1, 1}- e_{2, 6, 5} xf_{1, 1, 4}+ e_{3, 3, 3} xf_{4, 3, 3}+$$
$$ e_{3, 1, 2} xf_{4, 2, 1}- e_{2, 6, 2} xf_{1, 2, 4}- e_{2, 6, 5}
xf_{1, 2, 4}- e_{2, 2, 2} xf_{3, 4, 2}+ e_{4, 1, 1} xf_{1, 1, 1}-
e_{4, 1, 2} xf_{1, 1, 1}+ e_{1, 1, 2} xf_{1, 2, 1}- e_{4, 1, 2}
xf_{4, 3, 3}-$$ $$ e_{3, 4, 2} xf_{3, 2, 4}+ e_{2, 6, 2} xf_{1, 3,
4}+ e_{2, 6, 4} xf_{1, 3, 4}+ e_{2, 1, 1} xf_{1, 4, 4}- e_{1, 1, 2}
xf_{1, 4, 4}+ e_{1, 1, 1} xf_{1, 4, 4}- e_{2, 2, 4} xf_{3, 4, 2}+
e_{2, 3, 3} xf_{1, 4, 4}+$$ $$ e_{2, 6, 2} xf_{1, 4, 4}- e_{2, 1, 2}
xf_{1, 4, 4}- e_{2, 6, 1} xf_{1, 4, 4}- e_{2, 3, 5} xf_{1, 4, 4}-
e_{3, 1, 2} xf_{1, 4, 4}+ e_{3, 1, 1} xf_{1, 4, 4}+ e_{2, 3, 5}
xf_{1, 2, 1}+ e_{4, 1, 1} xf_{1, 4, 4}-$$ $$ e_{2, 1, 5} xf_{1, 2,
1}+ e_{3, 3, 3} xf_{1, 4, 4}+ e_{2, 3, 2} xf_{3, 2, 3}- e_{2, 1, 2}
xf_{2, 1, 1}- e_{1, 1, 2} xf_{2, 1, 1}- e_{4, 1, 2} xf_{1, 4, 4}+
e_{2, 5, 5} xf_{2, 1, 1}- e_{2, 3, 5} xf_{2, 1, 1}-$$ $$ e_{2, 6, 4}
xf_{2, 1, 1}- e_{3, 4, 3} xf_{3, 3, 4}- e_{1, 1, 2} xf_{3, 4, 4}-
e_{3, 1, 2} xf_{2, 1, 1}+ e_{2, 3, 3} xf_{3, 4, 4}+ e_{2, 1, 1}
xf_{2, 2, 1}+ e_{1, 1, 1} xf_{2, 2, 1}- e_{4, 1, 2} xf_{2, 1, 1}-$$
$$ e_{2, 1, 2} xf_{3, 4, 4}+ e_{2, 5, 5} xf_{3, 4, 4}+ e_{2, 4, 4}
xf_{2, 2, 1}+ e_{2, 3, 3} xf_{2, 2, 1}- e_{2, 6, 4} xf_{3, 4, 4}-
e_{3, 1, 2} xf_{3, 4, 4}+ e_{3, 1, 1} xf_{2, 2, 1}+ e_{2, 6, 6}
xf_{2, 2, 1}+$$ $$ e_{2, 5, 5} xf_{2, 2, 1}+ e_{3, 3, 3} xf_{2, 2,
1}+ e_{3, 4, 4} xf_{3, 4, 4}+ e_{2, 3, 5} xf_{4, 1, 3}+ e_{3, 3, 3}
xf_{3, 4, 4}+ e_{3, 4, 2} xf_{3, 4, 4}+ e_{1, 2, 1} xf_{2, 3, 1}-
e_{4, 1, 2} xf_{3, 4, 4}+$$ $$ e_{1, 1, 1} xf_{2, 3, 1}+ e_{4, 1, 1}
xf_{2, 2, 1}+ e_{3, 4, 4} xf_{2, 2, 1}+ e_{2, 5, 5} xf_{2, 3, 1}+
e_{2, 4, 4} xf_{2, 3, 1}+ e_{2, 3, 3} xf_{2, 3, 1}- e_{1, 1, 2}
xf_{1, 3, 1}+ e_{2, 2, 1} xf_{2, 3, 1}+$$ $$ e_{2, 1, 1} xf_{2, 3,
1}+ e_{3, 3, 3} xf_{2, 3, 1}+ e_{3, 2, 1} xf_{2, 3, 1}+ e_{3, 1, 1}
xf_{2, 3, 1}+ e_{2, 6, 6} xf_{2, 3, 1}+ e_{2, 3, 5} xf_{3, 2, 3}+
e_{4, 2, 1} xf_{2, 3, 1}+ e_{4, 1, 1} xf_{2, 3, 1}+$$ $$ e_{3, 4, 4}
xf_{2, 3, 1}+ e_{2, 3, 5} xf_{2, 1, 2}- e_{2, 1, 2} xf_{2, 2, 2}-
e_{1, 1, 2} xf_{2, 2, 2}+ e_{2, 5, 2} xf_{3, 1, 5}- e_{2, 6, 4}
xf_{2, 2, 2}- e_{2, 5, 1} xf_{3, 1, 5}- e_{2, 5, 3} xf_{3, 1, 5}-$$
$$ e_{3, 1, 2} xf_{2, 2, 2}- e_{4, 1, 2} xf_{2, 2, 2}+ e_{2, 5, 5}
xf_{3, 1, 5}- e_{2, 3, 3} xf_{2, 3, 2}+ e_{2, 3, 5} xf_{2, 3, 2}-
e_{3, 3, 3} xf_{2, 3, 2}- e_{1, 1, 2} xf_{1, 4, 1}+ e_{2, 5, 2}
xf_{3, 2, 5}+ e_{2, 5, 5} xf_{3, 2, 5}-$$ $$ e_{3, 3, 1} xf_{4, 1,
3}+ e_{2, 6, 4} xf_{2, 2, 3}- e_{1, 1, 2} xf_{3, 3, 3}+ e_{2, 1, 1}
xf_{3, 3, 5}+ e_{1, 1, 1} xf_{3, 3, 5}+ e_{3, 3, 2} xf_{4, 1, 3}-
e_{2, 1, 2} xf_{2, 3, 3}+ e_{2, 4, 4} xf_{3, 3, 5}+$$ $$ e_{2, 1, 1}
xf_{2, 3, 3}- e_{1, 1, 2} xf_{2, 3, 3}+ e_{2, 5, 3} xf_{3, 3, 5}+
e_{1, 1, 1} xf_{2, 3, 3}+ e_{2, 3, 3} xf_{3, 3, 5}- e_{2, 1, 2}
xf_{3, 3, 3}+ e_{2, 3, 3} xf_{2, 3, 3}- e_{3, 1, 2} xf_{2, 3, 3}+$$
$$ e_{3, 1, 1} xf_{2, 3, 3}- e_{2, 3, 5} xf_{2, 3, 3}+ e_{3, 3, 3}
xf_{3, 3, 5}+ e_{4, 1, 1} xf_{2, 3, 3}+ e_{3, 3, 3} xf_{2, 3, 3}+
e_{4, 1, 1} xf_{3, 3, 5}+ e_{3, 4, 4} xf_{3, 3, 5}- e_{4, 1, 2}
xf_{2, 3, 3}+$$ $$ e_{1, 2, 2} xf_{1, 1, 2}+ e_{2, 2, 1} xf_{1, 1,
2}+ e_{2, 1, 1} xf_{1, 1, 2}- e_{2, 5, 2} xf_{3, 4, 5}- e_{3, 3, 3}
xf_{4, 1, 3}- e_{2, 5, 4} xf_{3, 4, 5}+ e_{3, 1, 1} xf_{3, 1, 1}+
e_{2, 3, 3} xf_{3, 3, 3}-$$ $$ e_{3, 1, 2} xf_{3, 1, 1}+ e_{2, 3, 3}
xf_{1, 1, 2}+ e_{2, 2, 3} xf_{1, 1, 2}+ e_{4, 1, 1} xf_{3, 1, 1}-
e_{4, 1, 2} xf_{3, 1, 1}+ e_{1, 1, 2} xf_{3, 2, 1}- e_{1, 1, 2}
xf_{3, 5, 5}+ e_{2, 4, 4} xf_{1, 1, 2}+$$ $$ e_{2, 1, 2} xf_{3, 2,
1}- e_{2, 2, 5} xf_{1, 1, 2}- e_{2, 1, 2} xf_{3, 5, 5}- e_{2, 6, 4}
xf_{3, 3, 3}+ e_{2, 5, 5} xf_{1, 1, 2}+ e_{2, 3, 5} xf_{3, 2, 1}-
e_{2, 3, 5} xf_{3, 3, 3}- e_{2, 3, 5} xf_{3, 5, 5}+$$ $$ e_{2, 6, 6}
xf_{1, 1, 2}+ e_{3, 1, 1} xf_{1, 1, 2}+ e_{3, 2, 1} xf_{1, 1, 2}-
e_{3, 1, 2} xf_{3, 5, 5}+ e_{3, 2, 3} xf_{1, 1, 2}+ e_{3, 3, 3}
xf_{1, 1, 2}- e_{4, 1, 2} xf_{3, 5, 5}- e_{3, 2, 3} xf_{3, 3, 2}+$$
$$ e_{3, 4, 4} xf_{1, 1, 2}+ e_{4, 1, 1} xf_{1, 1, 2}+ e_{4, 2, 1}
xf_{1, 1, 2}- e_{1, 1, 2} xf_{1, 2, 2}+ e_{1, 2, 2} xf_{1, 2, 2}-
e_{3, 1, 2} xf_{3, 3, 3}- e_{3, 1, 3} xf_{3, 3, 1}- e_{2, 1, 2}
xf_{1, 2, 2}-$$ $$ e_{1, 1, 2} xf_{3, 4, 1}- e_{2, 1, 2} xf_{3, 4,
1}- e_{2, 6, 4} xf_{1, 2, 2}+ e_{4, 1, 1} xf_{4, 1, 1}- e_{2, 2, 5}
xf_{1, 2, 2}- e_{2, 1, 4} xf_{3, 4, 1}- e_{4, 1, 2} xf_{4, 1, 1}-
e_{2, 3, 5} xf_{1, 2, 2}+$$ $$ e_{2, 5, 5} xf_{1, 2, 2}+ e_{1, 1, 2}
xf_{4, 2, 1}+ e_{2, 1, 2} xf_{4, 2, 1}- e_{4, 1, 2} xf_{3, 3, 3}-
e_{3, 1, 2} xf_{1, 2, 2}+ e_{2, 3, 5} xf_{4, 2, 1}- e_{4, 1, 2}
xf_{1, 2, 2}- e_{1, 2, 2} xf_{1, 3, 2}+$$ $$ e_{1, 1, 1} xf_{3, 1,
2}+ e_{2, 1, 1} xf_{3, 1, 2}+ e_{1, 2, 2} xf_{3, 1, 2}+ e_{1, 1, 1}
xf_{4, 1, 2}+ e_{2, 1, 1} xf_{4, 1, 2}+ e_{1, 2, 2} xf_{4, 1, 2}+
e_{2, 2, 2} xf_{3, 1, 2}+ e_{2, 2, 2} xf_{4, 1, 2}-$$ $$ e_{2, 3, 2}
xf_{3, 4, 3}+ e_{2, 3, 3} xf_{3, 1, 2}+ e_{2, 3, 3} xf_{4, 1, 2}+
e_{2, 4, 4} xf_{3, 1, 2}+ e_{2, 4, 4} xf_{4, 1, 2}+ e_{2, 5, 5}
xf_{4, 1, 2}+ e_{2, 5, 5} xf_{3, 1, 2}+ e_{3, 1, 1} xf_{4, 1, 2}+$$
$$ e_{2, 6, 6} xf_{4, 1, 2}+ e_{3, 2, 2} xf_{4, 1, 2}+ e_{3, 1, 1}
xf_{3, 1, 2}+ e_{2, 6, 6} xf_{3, 1, 2}+ e_{3, 2, 1} xf_{3, 1, 2}+
e_{3, 3, 3} xf_{4, 1, 2}+ e_{3, 4, 4} xf_{4, 1, 2}+ e_{4, 1, 1}
xf_{4, 1, 2}+$$ $$ e_{3, 2, 3} xf_{3, 1, 2}+ e_{4, 2, 1} xf_{4, 1,
2}+ e_{3, 3, 3} xf_{3, 1, 2}+ e_{4, 1, 1} xf_{3, 1, 2}+ e_{3, 4, 4}
xf_{3, 1, 2}+ e_{1, 2, 2} xf_{4, 2, 2}- e_{1, 1, 2} xf_{4, 2, 2}-
e_{2, 3, 4} xf_{3, 4, 3}+$$ $$ e_{4, 2, 1} xf_{3, 1, 2}+ e_{1, 2, 1}
xf_{1, 4, 2}- e_{2, 1, 2} xf_{4, 2, 2}+ e_{2, 2, 2} xf_{4, 2, 2}+
e_{2, 3, 2} xf_{4, 2, 3}- e_{1, 2, 2} xf_{1, 4, 2}+ e_{1, 2, 2}
xf_{3, 2, 2}- e_{1, 1, 2} xf_{3, 2, 2}-$$ $$ e_{2, 6, 4} xf_{4, 2,
2}+ e_{2, 5, 5} xf_{4, 2, 2}- e_{2, 3, 5} xf_{4, 2, 2}- e_{2, 1, 2}
xf_{3, 2, 2}+ e_{2, 2, 2} xf_{3, 2, 2}- e_{3, 1, 2} xf_{4, 2, 2}+
e_{3, 2, 2} xf_{4, 2, 2}- e_{1, 2, 2} xf_{3, 4, 2}-$$ $$ e_{2, 6, 4}
xf_{3, 2, 2}- e_{4, 1, 2} xf_{4, 2, 2}- e_{2, 3, 1} xf_{3, 1, 3}-
e_{2, 3, 5} xf_{3, 2, 2}+ e_{2, 5, 5} xf_{3, 2, 2}+ e_{2, 3, 2}
xf_{3, 1, 3}- e_{3, 1, 2} xf_{3, 2, 2}+ e_{2, 4, 1} xf_{1, 1, 3}-$$
$$ e_{4, 1, 2} xf_{3, 2, 2}- e_{2, 3, 3} xf_{3, 1, 3}- e_{2, 4, 2}
xf_{1, 1, 3}+ e_{2, 4, 3} xf_{1, 1, 3}+ e_{2, 3, 5} xf_{4, 2, 3}-
e_{2, 3, 1} xf_{4, 1, 3}- e_{2, 4, 5} xf_{1, 1, 3}+ e_{2, 3, 2}
xf_{4, 1, 3}-$$ $$ e_{2, 3, 3} xf_{4, 1, 3})+\lambda^2 e_{2, 6, 4}
xf_{3, 3, 5}+e_{1, 2, 1} xf_{3, 4, 2}-e_{2, 3, 3} xf_{2, 1, 2}-e_{3,
3, 3} xf_{2, 1, 2}+e_{1, 1, 1} xf_{2, 2, 2}+e_{2, 4, 4} xf_{4, 2,
2}+e_{2, 1, 1} xf_{2, 2, 2}+$$ $$e_{2, 6, 6} xf_{2, 2, 2}+e_{3, 1,
1} xf_{2, 2, 2}+e_{4, 1, 1} xf_{2, 2, 2}+e_{3, 4, 4} xf_{4, 2,
2}-e_{2, 3, 3} xf_{4, 2, 3}+e_{4, 2, 1} xf_{4, 2, 2}+e_{2, 4, 4}
xf_{3, 2, 2}+e_{2, 2, 1} xf_{3, 4, 2}+$$ $$e_{2, 6, 4} xf_{2, 1,
3}-e_{2, 6, 6} xf_{2, 1, 3}+e_{3, 2, 1} xf_{3, 2, 2}+e_{2, 3, 6}
xf_{3, 4, 3}+e_{3, 2, 3} xf_{3, 2, 2}-e_{2, 6, 6} xf_{2, 2, 3}+e_{3,
4, 4} xf_{3, 2, 2}+e_{4, 2, 1} xf_{3, 2, 2}-$$ $$e_{2, 3, 1} xf_{4,
2, 3}+e_{4, 2, 1} xf_{2, 2, 1}+e_{3, 2, 1} xf_{2, 2, 1}+e_{2, 4, 4}
xf_{2, 1, 1}+e_{3, 2, 1} xf_{2, 1, 1}-e_{3, 1, 3} xf_{3, 5, 1}+e_{3,
4, 4} xf_{2, 1, 1}+e_{4, 2, 1} xf_{2, 1, 1}+$$ $$e_{1, 2, 1} xf_{2,
2, 1}+e_{2, 2, 1} xf_{2, 2, 1}-e_{3, 2, 3} xf_{3, 5, 2}+e_{2, 2, 1}
xf_{2, 1, 1}+e_{1, 2, 1} xf_{2, 1, 1}-e_{2, 6, 6} xf_{4, 2, 1}-e_{3,
1, 1} xf_{4, 2, 1}-e_{3, 3, 3} xf_{4, 2, 1}+$$ $$e_{2, 3, 1} xf_{3,
4, 3}+e_{2, 6, 4} xf_{4, 2, 1}-e_{2, 6, 6} xf_{1, 3, 4}-e_{2, 3, 3}
xf_{4, 2, 1}-e_{2, 1, 1} xf_{4, 2, 1}+e_{2, 1, 6} xf_{3, 4, 1}-e_{2,
6, 1} xf_{1, 3, 4}-e_{1, 1, 1} xf_{4, 2, 1}+$$ $$e_{2, 6, 1} xf_{1,
2, 4}+e_{4, 1, 1} xf_{3, 3, 3}+e_{2, 6, 3} xf_{1, 2, 4}+e_{2, 1, 1}
xf_{3, 4, 1}+e_{1, 1, 1} xf_{3, 4, 1}+e_{2, 6, 6} xf_{3, 3, 3}+e_{4,
2, 1} xf_{1, 3, 3}+e_{3, 1, 1} xf_{3, 3, 3}-$$ $$e_{3, 4, 1} xf_{1,
3, 3}+e_{4, 2, 1} xf_{3, 5, 5}+e_{3, 2, 1} xf_{1, 3, 3}+e_{2, 2, 1}
xf_{1, 3, 3}-e_{2, 4, 1} xf_{1, 3, 3}+e_{3, 4, 4} xf_{3, 5, 5}+e_{2,
6, 6} xf_{3, 5, 5}+e_{3, 2, 1} xf_{3, 5, 5}+$$ $$e_{1, 2, 1} xf_{1,
3, 3}+e_{3, 4, 3} xf_{1, 2, 3}-e_{3, 3, 3} xf_{3, 2, 1}+e_{2, 6, 4}
xf_{3, 2, 1}+e_{2, 4, 1} xf_{1, 2, 3}+e_{2, 4, 3} xf_{1, 2, 3}-e_{2,
6, 6} xf_{3, 2, 1}+e_{2, 4, 4} xf_{3, 5, 5}+$$ $$e_{3, 1, 3} xf_{3,
2, 1}+e_{3, 4, 1} xf_{1, 2, 3}+e_{2, 5, 3} xf_{3, 5, 5}-e_{2, 3, 3}
xf_{3, 2, 1}+e_{2, 2, 1} xf_{3, 5, 5}-e_{2, 1, 1} xf_{3, 2, 1}+e_{1,
2, 1} xf_{3, 5, 5}-e_{1, 1, 1} xf_{3, 2, 1}+$$ $$e_{2, 5, 6} xf_{3,
4, 5}+e_{2, 5, 1} xf_{3, 4, 5}+e_{4, 2, 1} xf_{3, 3, 5}+e_{4, 2, 1}
xf_{1, 2, 2}+e_{1, 2, 1} xf_{1, 3, 2}+e_{3, 2, 1} xf_{3, 3, 5}+e_{3,
2, 1} xf_{1, 2, 2}+e_{3, 2, 3} xf_{1, 2, 2}+$$ $$e_{3, 4, 4} xf_{1,
2, 2}+e_{2, 2, 3} xf_{1, 2, 2}+e_{2, 4, 4} xf_{1, 2, 2}+e_{2, 2, 1}
xf_{1, 2, 2}+e_{1, 2, 1} xf_{3, 3, 5}+e_{2, 1, 1} xf_{3, 3, 3}+e_{2,
2, 1} xf_{3, 3, 5}+e_{1, 1, 1} xf_{3, 3, 3}-$$ $$e_{2, 5, 3} xf_{3,
2, 5}-e_{3, 3, 3} xf_{1, 2, 1}+e_{3, 1, 3} xf_{1, 2, 1}-e_{2, 5, 1}
xf_{3, 2, 5}+e_{1, 1, 1} xf_{1, 3, 1}-e_{1, 1, 1} xf_{1, 2, 1}+e_{2,
6, 4} xf_{1, 2, 1}-e_{2, 3, 3} xf_{1, 2, 1}+$$ $$e_{2, 1, 3} xf_{1,
2, 1}-e_{2, 6, 6} xf_{1, 2, 1}-e_{3, 4, 3} xf_{3, 5, 4}+e_{4, 2, 1}
xf_{3, 4, 4}+e_{2, 2, 6} xf_{3, 4, 2}-e_{3, 4, 1} xf_{3, 4, 4}+e_{3,
2, 1} xf_{3, 4, 4}+e_{2, 4, 6} xf_{3, 4, 4}-$$ $$e_{2, 3, 3} xf_{3,
2, 3}+e_{2, 2, 1} xf_{3, 4, 4}+e_{1, 2, 1} xf_{3, 4, 4}-e_{2, 3, 1}
xf_{3, 2, 3}+e_{3, 4, 1} xf_{3, 2, 4}+e_{3, 4, 3} xf_{3, 2, 4}+e_{4,
1, 1} xf_{4, 3, 3}+e_{3, 1, 1} xf_{4, 3, 3}+$$ $$e_{2, 6, 6} xf_{4,
3, 3}+e_{2, 1, 1} xf_{4, 3, 3}-e_{3, 3, 3} xf_{3, 5, 3}+e_{1, 1, 1}
xf_{4, 3, 3}-e_{3, 3, 3} xf_{4, 2, 3}-e_{3, 3, 1} xf_{4, 2,
3}+\lambda xf_5 e_{4, 1, 2}.
$$

For each generic value of $K$ the algebra $U(\cal L)$ has the
following irreducible representation $V$. There are two bases
$\{v_{1,i_1,j_1},v_{2,i_2,j_2},v_{3,i_3,j_3}, v_{4,i_4,j_4} ; 1\le
i_1,j_1,i_4,j_4\le 2, 1\le i_2,j_2\le 6, 1\le i_3,j_3\le 4\}$ and
$\{ w_{1,i_1,j_1},w_{2,i_2,j_2},w_{3,i_3,j_3},w_{4,i_4,j_4},w_5;
1\le i_1,j_1\le 4, 1\le i_2,j_2,i_4,j_4\le 3, 1\le i_3,j_3\le 5\}$
of the space $V$ such that
$$e_{\alpha,i,j}v_{\beta,i^{\prime},j^{\prime}}=\delta_{\alpha,\beta}\delta_{j,i^{\prime}}v_{\alpha,i,j^{\prime}},$$
$$f_5w_5=w_5, \quad f_5w_{\alpha,i,j}=f_{\alpha,i,j}w_5=0, \quad
f_{\alpha,i,j}w_{\beta,i^{\prime},j^{\prime}}=\delta_{\alpha,\beta}\delta_{j,i^{\prime}}w_{\alpha,i,j^{\prime}}$$
and
$$v_{1,i,j}=w_{1,i,j}, \quad
v_{3,i,j}=s_{i,3}s_{j,3}w_{4,i,j}+w_{3,i,j}, \quad
v_{4,i,j}=s_{i,1}s_{j,1}w_5+w_{4,i,j},$$
$$v_{2,i,j}=s_{i,5}s_{j,5}w_{3,i,j}+\sum_{1\le
i^{\prime},j^{\prime}\le
3}\phi_i^{i^{\prime}}(\lambda)\phi_j^{j^{\prime}}(t)w_{2,i^{\prime},j^{\prime}}+\sum_{1\le
i^{\prime},j^{\prime}\le
4}\psi_i^{i^{\prime}}\psi_j^{j^{\prime}}w_{1,i^{\prime},j^{\prime}}.$$
Here $s_{k,l}=1$ if $k\le l$ and $s_{k,l}=0$ otherwise,
$\phi_k^{k^{\prime}}(u)=1$ if
$(k^{\prime},k)\in\{(1,2),(1,4),(1,5),(2,1),(2,3),$
$(2,4),(3,1),(3,6) \}$, $\phi_5^3(u)=u$ and
$\phi_k^{k^{\prime}}(u)=0$ otherwise, $\psi_k^{k^{\prime}}=1$ if
$(k^{\prime},k)\in\{(1,1),(2,2),(3,4),(4,6) \}$ and
$\psi_k^{k^{\prime}}=0$ otherwise, $\lambda\in\C$ is a parameter of
the algebra $U(\cal L)$ and $t\in\C$ is a parameter of
representation. In this representation $K$ acts as $\displaystyle
\frac{\lambda(t-1)}{t-\lambda}$.

\vskip.3cm \noindent {\bf Acknowledgments.} The authors are grateful
to I.Z. Golubchik and M.A. Semenov-Tian-Shansky for useful
discussions. The research was supported by the Manchester Institute
for Mathematical Sciences (MIMS). The research was partially
supported by: RFBR grant 05-01-00189, NSh grants 1716.2003.1 and
2044.2003.2.

\end{document}